\documentclass[12pt]{article}
\usepackage{latexsym,amssymb,amsmath,a4wide,color}

\newtheorem{thm}{Theorem}
\newtheorem{lem}{Lemma}
\newtheorem{cor}{Corollary}
\newtheorem{conj}{Conjecture}
\newtheorem{prop}{Proposition}
\newtheorem{exer}{Exercise}
\newcommand{\ebox}{\hfill $\Box$\\\vspace{0.15cm}}
\newcommand{\pr}{{\bf Proof.}\ }
\newcommand{\bt}{\begin{thm}}
\newcommand{\et}{\end{thm}}
\newcommand{\bl}{\begin{lem}}
\newcommand{\el}{\end{lem}}
\newcommand{\bp}{\begin{prop}}
\newcommand{\ep}{\end{prop}}
\newcommand{\bc}{\begin{cor}}
\newcommand{\ec}{\end{cor}}
\newcommand{\bcj}{\begin{conj}}
\newcommand{\ecj}{\end{conj}}
\newcommand{\bex}{\begin{exer}}
\newcommand{\eex}{\end{exer}}
\newcommand{\bi}{\begin{itemize}}
\newcommand{\ei}{\end{itemize}}
\newcommand{\be}{\begin{equation}}
\newcommand{\ee}{\end{equation}}
\newcommand{\ben}{\begin{enumerate}}
\newcommand{\een}{\end{enumerate}}

\newcommand{\mt}{t\kern-0.035cm\char39\kern-0.03cm}
\newcommand{\ml}{l\kern-0.035cm\char39\kern-0.03cm}
\newcommand{\md}{d\kern-0.035cm\char39\kern-0.03cm}

\newcommand{\veps}{\varepsilon}

\newcommand{\PG}{{\rm P\Gamma L}}

\newcommand{\GF}{{\rm GF}}

\newcommand{\dia}{{\rm dia}}
\newcommand{\off}{{\rm off}}
\newcommand{\PSL}{{\rm PSL}}
\newcommand{\PGL}{{\rm PGL}}

\newcommand{\GL}{{\rm GL}}
\newcommand{\SL}{{\rm SL}}
\newcommand{\Aut}{{\rm Aut}}
\newcommand{\ovl}{\overline}
\newcommand{\noi}{\noindent}
\newcommand{\orb}{{\rm orb}}

\begin{document}

\title{\vspace{-2.3cm} Orientably-regular maps \\ on twisted linear fractional groups}
%\thanks{Partially supported by...}}

\author{}
\date{}
\maketitle

\begin{center}
\vspace{-1.3cm}

{\large Grahame Erskine} \\
\vspace{1.5mm} {\small Open University, Milton Keynes, U.K.}\\

\vspace{5mm}

{\large Katar\'ina Hri\v{n}\'akov\'a} \\
\vspace{1.5mm} {\small Slovak University of Technology, Bratislava, Slovakia}

\vspace{5mm}

{\large Jozef \v Sir\'a\v n} \\
\vspace{1.5mm} {\small
Open University, Milton Keynes, U.K., and \\ Slovak University of Technology, Bratislava, Slovakia}

\vspace{4mm}

\end{center}

\begin{abstract}
We present an enumeration of orientably-regular maps with automorphism group isomorphic to the twisted linear fractional group $M(q^2)$ for any odd prime power $q$.

\vskip 3mm

\noi {\em Keywords:} Orientably-regular map; Automorphism group; Twisted linear fractional group.

\end{abstract}

\vskip 3mm

%------------------------------------------------
\section{Introduction}
\smallskip

A (finite) orientably-regular map ${\cal M}$ is a cellular embedding of a connected graph in a compact, oriented surface, such that the group $\Aut^+({\cal M})$ of all orientation-preserving automorphisms of the embedding is transitive, and hence regular, on arcs of the embedded graph. In a regular map, all vertices have the same valency, say $k$, and all faces are bounded by closed walks of the same length, say, $\ell$; the map is then said to be of type $(k,\ell)$. The group ${\cal A}=\Aut^+({\cal M})$ is generated by two elements $x$ and $y$ of order $k$ and $\ell$ such that $x$ acts as a clockwise rotation of ${\cal M}$ about a vertex by $2\pi/k$ and $y$ acts as a clockwise rotation by $2\pi/\ell$ about the centre of a face incident with the vertex. The product $xy$ is then a rotation of ${\cal M}$ about the centre of an edge that is incident to both the vertex and the face.
\smallskip

Orientably-regular maps can be viewed as maps having the `highest level' of orientation-preserving symmetry among general maps (i.e., cellular embeddings of graphs). Regularity of ${\cal A}$ on the arc set of the embedded graph enables one to identify the map ${\cal M}$ with the triple $({\cal A},x,y)$ in such a way that arcs, edges, vertices and faces correspond to (say, left) cosets of the trivial group (that is, to elements of ${\cal A}$) and of the subgroups $\langle xy\rangle$, $\langle x\rangle$ and $\langle y\rangle$ of $ {\cal A}$. Incidence between arcs, edges, vertices and faces is given by non-empty intersection of the corresponding cosets, and the action of ${\cal A}$ on the cosets is simply given by left multiplication.
\smallskip

It follows that orientably-regular maps are, up to isomorphism, in a one-to-one correspondence with equivalence classes of triples $(G,x,y)$, where $G$ is a finite group admitting a presentation of the form $G=\langle x,y;\ x^k=y^\ell=(xy)^2=\ldots=1\rangle$, with two triples $(G_1,x_1,y_1)$ and $(G_2,x_2,y_2)$ being equivalent if there is a group isomorphism $G_1\to G_2$ taking $x_1$ onto $x_2$ and $y_1$ onto $y_2$. This way, investigation of orientably-regular maps can be reduced to purely group-theoretic considerations. The corresponding algebraic theory has been developed in depth in the influential paper \cite{JoSi1}.
\smallskip

The study of orientably-regular maps has rich history and, save Platonic solids and Kepler's polyhedra, has roots in the late 19th century; for historical information and surveys see e.g. \cite{JoSi1,JoSi2,Ne,Sir1}. For an excellent introduction into fascinating connections between orientably-regular maps, Dyck's triangle groups, Riemann surfaces and Galois groups we recommend \cite{JoSi2}. In particular, classification of orientably-regular maps has implications in classification of Riemann surfaces, cf. \cite{JoSi1,JoSi2}.
\smallskip

Since the concept of an orientably-regular map includes the underlying graph, the carrier surface and the supporting automorphism group, classification attempts for such maps in most cases follow one of these three directions. A number of influential results have been obtained in classification of orientably-regular maps in the first two directions; we refer to \cite{Sir1} for the most recent survey. In this paper we focus on the third direction, that is, classification of orientably-regular maps by their automorphism groups, in which results are much less abundant.
\smallskip

Leaving the trivial case of Abelian groups aside, classification of orientably-regular maps with a given automorphism group has been completed for the following types of groups, ranked by complexity of their structure:
\bi
\item groups of nilpotency class two and three \cite{Ban+};
\item products of two cyclic groups one of which acts regularly on vertices of the map  \cite{CT};
\item groups with cyclic odd-order and dihedral even-order Sylow subgroups \cite{CPS};
\item $\PSL(2,q)$ and $\PGL(2,q)$, where $q$ is an arbitrary prime power \cite{Mac,Sah,CPS2};
\item Ree groups $^2{\rm }G_2(3^n)$ for odd $n>1$, restricted to maps of type $(k,\ell)$ for $\ell=3$ and $k=7,9$ and all prime $k\equiv 11$ mod $12$ \cite{Jo-Ree};
\item Suzuki group ${\rm Sz}(2^n)$ for odd $n>1$, restricted to maps of type $(5,4)$ \cite{Jo-Suz}.
\ei
By this list, the only family of {\em simple} groups for which a classification of the corresponding orientably-regular maps is known are the groups $\PSL(2,q)$ for any prime power $q>3$. Classification of orientably-regular maps with automorphism group $\PGL(2,q)$, the obvious degree-two extension of $\PSL(2,q)$, can be extracted from the corresponding classification for $\PSL(2,q^2)$ through the well-understood inclusion $\PGL(2,q)<\PSL(2,q^2)$, cf. \cite{Sah,CPS2}.
\smallskip

For odd $q$, however, the simple group $\PSL(2,q^2)$ admits another interesting extension of degree two, namely, the group $M(q^2)$, also known as a {\em twisted} linear fractional group. By a classical result of Zassenhaus \cite{Zas}, the groups $\PGL(2,q)$ and $M(q^2)$ are the only finite sharply $3$-transitive groups (of degree $q+1$ and $q^2+1$, respectively). This motivates the question of classification of orientably-regular maps with automorphism group isomorphic to $M(q^2)$.
\smallskip

In this paper we present a complete enumeration of (isomorphism classes of) orientably-regular maps with automorphism group isomorphic to $M(q^2)$. The results are strikingly different from those for the groups $\PGL(2,q)$ in many ways. To give three examples, we note that (a) all the orientably-regular maps for $\PGL(2,q)$ are reflexible, while this is not the case for $M(q^2)$; (b) the groups $\PGL(2,q)$ are also automorphism groups of non-orientable regular maps while the groups $M(q^2)$ are not; and (c) for any even $k,\ell\ge 4$ not both equal to $4$ there are orientably-regular maps of type $(k,\ell)$ with automorphism group $\PGL(2,q)$ for infinitely many values of $q$, while for infinitely many such pairs $(k,\ell)$ there are no orientably-regular maps for $M(q^2)$ of that type for any $q$.
\smallskip

By our outline and the algebraic theory of \cite{JoSi1}, enumeration of orientably-regular maps with a given automorphism group $G$ reduces to enumeration of all triples $(G,x,y)$ with $G=\langle x,y;\ x^k=y^\ell=(xy)^2=\ldots=1\rangle$ up to conjugation by elements of $\Aut(G)$, that is, by considering triples $(G,x,y)$ and $(G,x',y')$ equivalent if there is an automorphism of $G$ taking $(x,y)$ onto $(x',y')$. We do this systematically for the twisted linear groups $G=M(q^2)$. In sections \ref{sec:pre} and \ref{s:tw-sbgp} we introduce the group $M(q^2)$ and study its subgroups. Sections \ref{sec:diag}, \ref{sec:twi} and \ref{sec:twist} deal with identifying `canonical' forms of elements of $G$ and study their conjugacy in depth. In sections \ref{sec:twpairs}, \ref{sec:rep-dia} and \ref{sec:rep-off} we develop arguments for counting `canonical' pairs of elements of $G$. All the auxiliary facts are then processed in section \ref{sec:enum} to produce our main result:
\bigskip

\noi {\bf Theorem} {\sl Let $q=p^f$ be an odd prime power, with $f=2^{\alpha}o$ where $o$ is odd. The number of orientably-regular maps $\cal M$ with $\Aut^+({\cal M})\cong M(q^2)$ is, up to isomorphism, equal to      \[ \frac{1}{f}\sum_{d\mid o} \mu(o{/}d)h(2^{\alpha}{}d)\ ,\] where $h(x)= (p^{2x}-1)(p^{2x}-2)/8$ and $\mu$ is the M\"obius function.}
\bigskip

\noi We note that this result may be interpreted as counting generating pairs $(x,y)$ of $G=M(q^2)$ such that $(xy)^2=1$, up to conjugacy in $\Aut(G)$, which may be of independent interest to specialists in group theory.
The final section \ref{sec:rem} contains related results and remarks.
\smallskip

%------------------------------------------------
\section{The twisted linear groups {\bf $M(q^2)$}}\label{sec:pre}
\smallskip

For a finite field $F$ let $S(F)$ and $N(F)$ be the set of non-zero squares and non-squares of $F$. The general linear group $\GL(2,F)$ is the group of all non-singular $2\times 2$ matrices with entries in $F$; restriction to matrices with determinant $1$ gives the special linear group $\SL(2,F)$. The groups $\PGL(2,F)$ and $\PSL(2,F)$, the quotients of $\GL(2,F)$ and $\SL(2,F)$ by the corresponding centres, are known as the linear fractional groups. They can equivalently be described as groups of all transformations $z\mapsto (az+b)/(cz+d)$ of the set $F\cup\{\infty\}$ (with the obvious rules for calculations with $\infty$), with $ad-bc\ne 0$ and $ad-bc\in S(F)$ for $\PGL(2,F)$ and $\PSL(2,F)$, respectively. The group $\PSL(2,F)$ is an index $2$ subgroup of $\PGL(2,F)$ unless $F$ has characteristic $2$, in which case the two groups are the same.
\smallskip

Suppose now that $F$ admits an automorphism $\sigma$ of order $2$, which happens if and only if $|F|=q^2$ for some prime power $q$, and $\sigma$ is then given by $x\mapsto x^q$ for every $x\in F$. If, in addition, $q$ is odd, then one may `twist' the transformations described above by considering the permutations of $F\cup\{\infty\}$ defined by $z\mapsto (az+b)/(zc+d)$ if $ad-bc\in S(F)$ and $z\mapsto (az^{\sigma}+b)/(cz^{\sigma}+d)$ if $ad-bc\in N(F)$. These transformations form a group under composition, denoted $M(F)$ or $M(q^2)$, and called the {\em twisted fractional linear group}. Observe that $\PSL(2,F)$ is a subgroup of $M(F)$ of index two, again. By a well-known result due to Zassenhaus \cite{Zas}, the groups $\PGL(2,F)$ for an arbitrary finite field $F$, and $M(F)$ for fields of order $q^2$ for an odd prime power $q$, are precisely the finite, sharply $3$-transitive permutation groups (on the set $F\cup\{\infty\}$ in both cases).
\smallskip

In this paper we will focus on the twisted fractional linear groups, with the goal to classify the orientably-regular maps they support. For our purposes, however, it will be useful to work with a different representation of these groups. From this point on, let $F=\GF(q^2)$ for some odd prime power $q$ and let $F_0\cong \GF(q)$ be its unique subfield of order $q$; let $F^*$ and $F_0^*$ be the corresponding multiplicative groups. Further, let $\sigma$ be the unique automorphism of $F$ of order $2$; we have $x^{\sigma}=x^q$ for any $x\in F$, and $x^{\sigma}=x$ if and only if $x\in F_0$. If $A\in \GL(2,F)$, by $A^{\sigma}$ we denote the matrix in $\GL(2,F)$ obtained by applying $\sigma$ to every entry of $A$.
\smallskip

Let $J=\GL(2,F)\rtimes Z_2$, where multiplication in the semidirect product is defined by
$(A,i)(B,j)=(AB^{\sigma^{i}},i+j)$; equivalently, $J$ is an extension of $\GL(2,F)$ by the automorphism $\sigma$. To introduce a `twisted' subgroup of $J$, for every $A\in \GL(2,F)$ we define $\iota_A\in Z_2=\{0,1\}$ by letting $\iota_A=0$ if $\det(A)\in S(F)$ and $\iota_A=1$ if $\det(A)\in N(F)$. We now let $K=\{(A,\iota_A);\ A\in \GL(2,F)\}$; multiplication in $K$ is, of course, given by $(A,\iota_A)(B,\iota_B) =(AB^{\sigma^{\iota_A}},\iota_A+\iota_B)$ for any $A,B\in \GL(2,F)$. The group $K$ and its quotient groups will be of principal importance in what follows.
\smallskip

Let $K_0=\{(A,0);\ A\in \GL(2,F),\ \iota_A=0\}$ be the subgroup of $K$ index $2$ of $K$. The centre $L$ of $K_0$ consists of elements of the form $(D,0)$, where $D\in \GL(F)$ is a scalar matrix; obviously $L$ is also a normal subgroup of both $K$ and $J$. It can be checked that the factor group $G=K/L$ is isomorphic to $M(q^2)$, and since $K$ has index $2$ in $J$, the group $G=M(q^2)$ is (isomorphic to) a subgroup of index $2$ of $\ovl{G}=J/L$. The group $\ovl{G}$ can alternatively be described as $G\langle\sigma\rangle$, the split extension of $G$ by $\langle\sigma\rangle\simeq Z_2$. Observe also that the factor group $G_0=K_0/L$ is isomorphic to $\PSL(2,F)$, and if $q$ is a prime, the group $J/L$ is isomorphic to $\PG (2,q^2)$.
\smallskip

Elements $(A,i)L$, that is, cosets $\{(\delta A,i);\ \delta\in F^*\}$, of the factor groups $G=K/L$ and $\ovl{G}=J/L$ will throughout be denoted $[A,i]$; they will be called {\em untwisted} if $i=0$ and {\em twisted} if $i=1$.
\smallskip

For our final enumeration it will be necessary to determine the automorphism group of $M(q^2)$. While the result appears to be `obvious' we provide a simple proof based on an fact which may be folklore to group-theorists.
\smallskip

\bl\label{autgrp}
Let $U$ be a characteristic subgroup of a group $\tilde U$ of index $2$. Suppose that the centre of $U$ is trivial and every automorphism of $U$ extends to an automorphism of $\tilde U$. Then $\Aut(U)\cong Aut(\tilde U)$.
\el

\pr
The assumption of $U$ being characteristic in $\tilde U$ implies that every $h\in \Aut(\tilde U)$ restricts to an $h_U\in \Aut(U)$. Since each automorphism of $U$ extends to an automorphism of $\tilde U$, the assignment $h\mapsto h_U$ is a group epimorphism $\vartheta:\ \Aut(\tilde U)\to \Aut(U)$. Suppose that $h\in \Aut(\tilde U)$ is in the kernel of $\vartheta$, so that $h_U$ is the identity mapping on $U$. For every $x\in U$ and every $y\in \tilde U{\setminus}U$ we have $y^{-1}xy\in U$ and hence $y^{-1}xy=h(y^{-1}xy)=h(y)^{-1}xh(y)$, which implies that $h(y)y^{-1}$ commutes with $x$ for all $x\in U$. Observe that $h(y)y^{-1}\in U$, since $U$ was assumed to have index $2$ in $\tilde U$. By triviality of the centre of $U$ we have $h(y)y^{-1}=1$ and as this is valid for all $y\in \tilde U{\setminus}U$ we conclude that $h$ is the identity on $\tilde U$. It follows that the kernel of $\vartheta$ is trivial and so $\Aut(U)\cong Aut(\tilde U)$.
\ebox

We now apply Lemma \ref{autgrp} to $U=G_0={\rm PSL}(2,q^2)$ and $\tilde U=G=M(q^2)$ for $q=p^f$, where $p$ is an odd prime and $f$ a positive integer. Being a simple subgroup of $M(q^2)$, the group $G_0$ is characteristic (and of index two) in $G$. It is well known (see e.g. \cite{Giu}) that $\Aut(G_0)\cong {\rm P\Gamma L}(2,q^2) \simeq {\rm PGL}(2,q^2)\rtimes Z_{2f}$, with an element $(C,\varphi)\in {\rm PGL}(2,q^2)\rtimes Z_{2f}$ acting on $G_0$ by $X\mapsto (C^{-1}XC)^{\varphi}$. Now, any $(C,\varphi)$ is easily seen to extend to $G$ by $[X,\iota_X]\mapsto ([C,0]^{-1}[X,\iota_X][C,0])^{\varphi}$. By Lemma \ref{autgrp} we now obtain:

\bp\label{cor-auto}
The automorphism group of $M(q^2)$ is isomorphic to ${\rm P\Gamma L}(2,q^2)$. \hfill $\Box$
\ep

%------------------------------------------------
\section{Twisted subgroups of {\bf $M(q^2)$}}\label{s:tw-sbgp}
\smallskip

Let $q=p^f$ for an odd prime $p$ and a positive integer $f$; these will be fixed throughout. In this section we will focus on the {\em twisted} subgroups of $M(p^{2f})$, that is, those isomorphic to $M(p^{2e})$ for suitable $e\le f$. From now on we will use the notation $F_{m}=\GF(p^{m})$ for a Galois field of order $p^m$ for $m\le f$ but keep letting $F=\GF(p^{2f})$. We begin by identifying the possible values of $e$.

\bl\label{subgps}
A group $M(p^{2e})$ is isomorphic to a subgroup of $M(p^{2f})$ if and only if $e$ is a divisor of $f$ such that  $f/e$ is odd.
\el

\pr
If $e$ divides $f$, then $p^{2f}-1=(p^{2e}-1)Q$ for $Q=p^{2e(d-1)} + p^{2e(d-2)} + \ldots + p^{2e}+1$ and $d=f/e$. Assume that $d$ is odd. If $\xi$ is a primitive element of $F$, then $\xi^Q$ is a primitive element of $F_{2e}<F$. Further, the restriction of the assignment $x\mapsto x^{p^f}$ for $x\in F$ onto $F_{2e}$ coincides with the mapping $x\mapsto x^{p^e}$ for $x\in F_{2e}$. To see this, it is clearly sufficient to consider the effect on taking the powers of $x=\xi^Q$, for which we have $x^{p^e}=x^{p^f}$ if and only if $Qp^e\equiv Qp^f$ mod $(p^{2f}-1)$. The congruence is further equivalent to $p^e \equiv p^f$ mod $(p^{2e}-1)$ and also to $p^{2e}-1$ dividing $p^{f-e}-1$, which is true if and only if $2e$ is a divisor of $f-e=(d-1)e$, and the last condition is satisfied because $d$ is odd. This shows that the restriction of the automorphism $\sigma:\ x\mapsto x^{p^f}$ of $F$ onto $F_{2e}$ coincides with the automorphism $\sigma':\ x\mapsto x^{p^e}$ of $F_{2e}$. By construction of the twisted linear groups introduced in section \ref{sec:pre} it is now clear that $M(p^{2f})$ contains a copy of $M(p^{2e})$.
\smallskip

To prove the reverse implication, assume that $H\cong M(p^{2e})$ is a subgroup of $G=M(p^{2f})$. From the order of $H$ dividing the order of $G$ one sees that $e$ divides $f$, and our goal is to show that $f/e$ is odd. We prove this by induction on $f/e$. The case when $f/e=1$ is obvious and so we assume that $f>e$. Let $\tilde H$ be a maximal subgroup of $G$ containing $H$. From Theorem 1.5 of \cite{Giu} it follows that $\tilde H$ is isomorphic to the normaliser of ${\rm PSL}(2,p^{2g})$ in $G$ for some positive divisor $g$ of $f$ such that $f/g$ is an odd prime. As before, let $G_0$ be the (normal) subgroup of index $2$ in $G$ isomorphic to ${\rm PSL}(2,p^{2f})$. The group $\tilde H \cap G_0$ is contained in a maximal subgroup of $G_0$ containing ${\rm PSL}(2,p^{2g})$. However, the group ${\rm PSL}(2,p^{2g})$ itself is a maximal subgroup of $G$ as $f/g$ is odd, cf. \cite{Giu} again. It follows that $\tilde H \cap G_0$ must be isomorphic to ${\rm PSL}(2,p^{2g})$, and therefore $\tilde H$ must be isomorphic to $M(p^{2g})$. But $\tilde H$ contains also $M(p^{2e})\cong H$ and therefore $e$ divides $g$. As $g/e < f/e$, we may apply our induction hypothesis by which $g/e$ is odd. Since $f/g$ was an odd prime we conclude that $f/e$ is odd as well, completing the induction step.
\ebox

If $f/e$ is odd, a particularly important copy of $M(p^{2e})$ in $M(p^{2f})$ is formed by all the pairs $[X,\iota_X]$ with $X\in {\rm GL}(2,p^{2e})$ such that all entries of $X$ lie in the subfield $F_{2e}$ of $F$; this copy will be called {\em canonical}. The copy of ${\rm PSL}(2,p^{2e})$ in $M(p^{2f})$ formed by all the pairs $[X,0]$ with $X\in {\rm SL}(2,F_{2e})$ will be called {\em canonical} as well. We now prove a useful auxiliary result on canonical subgroups.

\bp\label{e/f-conjug} Let $f/e$ be an odd integer and let $H\cong M(p^{2e})$ be a subgroup of $G=M(p^{2f})$ such that $H$ contains the canonical copy of ${\rm PSL}(2,p^{2e})$. Then, $H$ is equal to the canonical copy of $M(p^{2e})$ in $G$.
\ep

\pr Let $H$ be a copy of $M(p^{2e})$ in $G$ such that $H_0 = H\cap {\rm PSL}(2,p^{2f})$ is equal to the canonical copy of ${\rm PSL}(2,p^{2e})$ in $G$. Obviously, $H_0$ is a normal subgroup of $H$ of index two. Let $[A,1]$ be an element of $H{\setminus}H_0$, where $A$ is the $2\times 2$ matrix with rows $(a,b)$ and $(c,d)$ for some $a,b,c,d\in F$ with $\delta = ad-bc\in N(F)$. We may assume that the entry $c$ in the lower left corner of $A$ is non-zero. Indeed, if $c=0$ and $b\ne 0$, letting $D$ be an off-diagonal matrix with entries $-1$ and $1$ we may replace $[A,1]$ with the product $[D,0][A,1]\in H{\setminus}H_0$, and if $A$ is a diagonal matrix we may replace $[A,1]$ with the product $[D',0][A,1]\in H{\setminus}H_0$ for a matrix $D'$ with rows $(1,0)$ and $(1,1)$. Then, since we are working with projective groups, we may assume that $c=1$, so that $\delta=ad-b$.
\smallskip

By our assumption the group $H_0$ also contains the element $[C,0]$ with $C$ having rows $(1,1)$ and $(0,1)$. Normality of $H_0$ in $H$ implies that $[A,1][C,0][A,1]^{-1}=[ACA^{-1},0]\in H_0$ and also $[A,1]^{-1}[C,0][A,1]=[(A^{\sigma})^{-1}CA^{\sigma},0]\in H_0$. Evaluating the products we obtain
\[ \veps ACA^{-1} = \begin{pmatrix} a & b \\ 1 & d\end{pmatrix}
\begin{pmatrix} 1 & 1 \\ 0 & 1 \end{pmatrix} \begin{pmatrix} d & -b \\ -1 & a\end{pmatrix} = \begin{pmatrix} \delta - a & a^2 \\ -1 & \delta +a \end{pmatrix}\ , \ {\rm and}  \]
\[ \veps'(A^{\sigma})^{-1}CA^{\sigma} = \begin{pmatrix} d^{\sigma} & -b^{\sigma} \\ -1 & a^{\sigma}\end{pmatrix} \begin{pmatrix} 1 & 1 \\ 0 & 1 \end{pmatrix} \begin{pmatrix} a^{\sigma} & b^{\sigma} \\ 1 & d^{\sigma}\end{pmatrix} =  \begin{pmatrix} \delta^{\sigma} + d^{\sigma} & (d^{\sigma})^2 \\ -1 & \delta^{\sigma} - d^{\sigma} \end{pmatrix}\  \]
for some $\veps,\veps'\in F^*$. Since $H_0$ is assumed to be equal to the canonical copy ${\rm PSL}(2,p^{2e})$  in $G$, all the remaining entries of the two matrices on the right-hand sides above must lie in $F_{2e}$. This readily implies that both $a, \delta, d^{\sigma} \in F_{2e}$, and since $F_{2e}$ is setwise preserved by $\sigma$ we also have $d\in F_{2e}$ and so $b=ad-\delta\in F_{2e}$ as well. We conclude that $A\in {\rm GL}(2,p^{2e})$ and hence the subgroup $H$ generated by $[A,1]$ and $H_0$ is identical with the canonical copy of $M(p^{2e})$ in $G$. \ebox

As a consequence we prove that all twisted subgroups of $G$ are conjugate. Recall that $G_0$ denotes the (unique) subgroup of $G$ isomorphic to ${\rm PSL}(2,p^{2f})$.

\bp\label{e/f-conjug-c}
If $f/e$ is an odd integer, then all subgroups of $G=M(p^{2f})$ isomorphic to $M(p^{2e})$ are conjugate in $G_0$.
\ep

\pr
Let $H\cong M(p^{2e})$ be a subgroup of $G$; letting $H_0=H\cap G$ we obviously have $H_0\cong {\rm PSL}(2,p^{2e})$. The classification of subgroups of ${\rm PSL}(2,p^{2f})$, nicely displayed in \cite{LS}, tells us that in the case when $f/e$ is odd, all subgroups of $G_0$ isomorphic to ${\rm PSL}(2,p^{2e})$ are conjugate in $G_0$. It follows that there is an inner automorphism $\pi$ of $G_0$ such that $\pi(H_0)$ is equal to the canonical copy of ${\rm PSL}(2,p^{2e})$ in $G_0$. The group $\pi(H)$ is then isomorphic to $M(p^{2e})$ and contains the canonical copy of ${\rm PSL}(2,p^{2e})$. By Proposition \ref{e/f-conjug}, however, the group $\pi(H)$ must be {\em equal} to the canonical copy of $M(p^{2e})$ in $G$. In particular, all subgroups of $G$ isomorphic to $M(p^{2e})$ are conjugate in $G_0$. \ebox

We conclude with a sufficient condition for a subgroup of $G$ to be twisted; this result will be of key importance later. In order to state it, we will say that a subgroup $H$ of $G$ {\em stabilises a point} if, in the natural action of $G$ on the set $F\cup \{\infty\}$ via linear fractional mappings from section \ref{sec:pre}, there exists a point in $F\cup \{\infty\}$ fixed by all linear fractional mappings corresponding to elements of $H$. Also, for any positive divisor $g$ of $2f$ let $G_{g}$ be the canonical copy of ${\rm PSL}(2,p^{g})$ in the group $G=M(p^{2f})$. Moreover, if $g$ is even, we let $G^*_{g}$ denote the copy of ${\rm PGL}(2,p^{g/2})$ in $G_{g}$ formed by (equivalence classes of) non-singular $2\times 2$ matrices over $\GF(p^{g/2})$.

\bp\label{subgps-class}
Let $H$ be a subgroup of $G=M(p^{2f})$ not contained in the subgroup $G_0={\rm PSL}(2,p^{2f})$ and let $H_0=H\cap G_0$. If $H_0$ does not stabilise a point and is neither dihedral nor isomorphic to $A_4$, $S_4$ or $A_5$, then $H$ is conjugate in $\ovl{G}$ to a subgroup isomorphic to $M(p^{2e})$ for some positive divisor $e$ of $f$ such that $f/e$ is odd.
\ep

\pr
By a handy summary \cite{LS} of Dickson's classification of subgroups of projective special linear groups over finite fields, subgroups of $G_0$ comprise point stabilisers, dihedral groups, $A_4$, $S_4$, $A_5$, and ${\rm PSL}(2,p^{g})$ for divisors $g$ of $2f$ together with ${\rm PGL}(2,p^{g/2})$ for even divisors $g$ of $2f$.
Our assumptions imply that $H_0$ must be isomorphic to one of the last two types of subgroups.
\smallskip

Let $e$ be the smallest positive integer such that $f/e$ is odd and $g$ divides $2e$. Invoking the classification of subgroups of $G_0$ again \cite{LS}, and specifically the fact that any two copies of ${\rm PSL}(2,p^{2e})$ in $G$ are conjugate (by an element of $G$), we may assume that $H_0$ is a subgroup of the canonical copy $G_{2e}$ of ${\rm PSL}(2,p^{2e})$ in $G$. Now, by \cite{LS}, the number of conjugacy classes of copies of ${\rm PSL}(2,p^{g})$ in $G_{2e}$ is $1$ or $2$ according to whether $2e/g$ is odd or even. The two classes if $g$ divides $e$ are, however, by \cite{Ca+}, fused under conjugacy in ${\rm PGL}(2,p^{2e})$ and hence also under conjugacy in $\ovl{G}$; the same holds for copies of ${\rm PGL}(2,p^{g/2})$ if $g$ is even. We therefore may assume that $H_0$ is {\em equal} either to the canonical copy $G_{g}$ of ${\rm PSL}(2,p^{g})$ in $G$, or to the canonical copy $G^*_{g}$ of ${\rm PGL}(2,p^{g/2})$ in $G$ if $g$ is even.
\smallskip

Let $[A,1]\in H{\setminus}H_0$ and $[C,0]\in H_0$ be the elements from the proof of Proposition \ref{e/f-conjug}. Considering the conjugation $[A,1][C,0][A,1]^{-1}=[ACA^{-1},0]\in H_0$ and verbatim repeating the arguments from the proof of Proposition \ref{e/f-conjug}, except with $F_{g}$ in place of $F_{2e}$, we conclude that $\delta\in F_{g}$. However, a non-square element $\delta\in F=F_{2f}$ can live in a subfield $F_{g} < F_{2f}$ only if $2f/g$ is odd. As $g\mid 2e$ with $f/e$ odd and $e$ was the smallest positive integer with this property, oddness of $2f/g$ implies that $g=2e$.
\smallskip

We have arrived at the conclusion that $H_0$ is {\em equal} either to the canonical copy $G_{2e}$ of ${\rm PSL}(2,p^{2e})$ or to the canonical copy $G^*_{2e}$ of ${\rm PGL}(2,p^{e})$. But in the second case, by the argument in the preceding paragraph, we would have $\delta\in F_e$ while $2f/e$ is not odd, a contradiction. Therefore $H_0$ must be equal to $G_{2e}$. The conclusion that $H\cong M(p^{2e})$ now follows from Proposition \ref{e/f-conjug}.
\ebox

%------------------------------------------------
\section{Representatives of twisted elements}\label{sec:diag}
\smallskip

Recalling the notation introduced in section \ref{sec:pre}, we begin by identifying elements in conjugacy classes of $K{\setminus}K_0$ that have a particularly simple form. To facilitate the description here and also in the sections that follow, we let $\dia(\alpha,\beta)$ and $\off(\alpha,\beta)$, respectively, denote the $2\times 2$ matrix with diagonal entries $\alpha,\beta$ (from the top left corner) and zero off-diagonal entries, and the  $2\times 2$ matrix with off-diagonal entries $\alpha,\beta$ (from the top right corner) and zero diagonal entries.
\smallskip

For every element $(A,1)\in K{\setminus}K_0$ we have $(A,1)^2=(AA^{\sigma},0)$. In the study of conjugacy in $K{\setminus}K_0$ it turns out to be  important to understand the behaviour of the products $AA^{\sigma}$. Observe that if $\delta=\det(A)\in N(F)$, then $\det(AA^{\sigma})=\delta\delta^{\sigma}\in N(F_0)$.
\smallskip

Let $(A,1)\in K{\setminus}K_0$ and let $\{\lambda_1,\lambda_2\}$ be the spectrum of $AA^{\sigma}$ in a smallest extension $F'$ of $F$ of degree at most two in which $\sigma$ may still be assumed to be given by $x\mapsto x^q$. Since $A^{\sigma}A$ is both a conjugate and also a $\sigma$-image of $AA^{\sigma}$, we have $\{\lambda_1,\lambda_2\}^{\sigma}= \{\lambda_1,\lambda_2\}$. This means that either (1) $\lambda_i^{\sigma}=\lambda_i$ for $i=1,2$, or (2) $\lambda_1^{\sigma}=\lambda_2$ and $\lambda_2^{\sigma} =\lambda_1$. Note that (2) implies $\lambda_1^{q^2}= (\lambda_1^q)^q=\lambda_2^q=\lambda_1$ and, similarly, $\lambda_2^{q^2} =\lambda_2$. We conclude that $F'=F$ in both the situations (1) and (2) and so both $\lambda_1,\lambda_2$ are in $F$. Observe that $\lambda_1\ne \lambda_2$, as otherwise we would have $\lambda_1=\lambda_1^{\sigma}\in F_0$ and $\det(AA^{\sigma})=\lambda_1^2\in S(F_0)$, a contradiction. Moreover, it follows that in the case (1) we have $\lambda_i\in F_0$ for $i=1,2$ with $\lambda_1\lambda_2\in N(F_0)$, and in the case (2) $\lambda_i\in N(F)\subset F{\setminus}F_0$ since $\det(AA^{\sigma})=\lambda_1\lambda_1^{\sigma}\in N(F_0)$.
\smallskip

We will now refine our considerations of $AA^{\sigma}$. As before, let $q=p^{f}$ for some odd prime $p$ and let $e$ be the {\em smallest} positive divisor of $f$ with $f/e$ odd such that $AA^{\sigma}=\veps C$ for some $C\in {\rm SL}(2,p^{2e})$ and for some $\veps\in F^*$. In other words, we look for the smallest subfield $F_{2e}$ of $F$, with $f/e$ odd, such that all entries of $C$ lie in $F_{2e}$; note that we may assume $C$ to have determinant $1$ since the determinant of $AA^{\sigma}$ is a non-zero square of $F$. If $\{\mu,\mu^{-1}\}$ is the spectrum of $C$, we have, without loss of generality, $\lambda_1=\veps\mu$ and $\lambda_2=\veps\mu^{-1}$. Observe that since $\lambda_1,\lambda_2,\veps\in F$, we have $\mu,\mu^{-1}\in F$. Now, $\mu,\mu^{-1}$ are roots of a quadratic polynomial over $F_{2e}$ and therefore both belong to $F_{2e}$ or to a quadratic extension of $F_{2e}$. But as $f/e$ is odd, the field $F$ does not contain a quadratic extension of $F_{2e}$. We conclude that $\mu,\mu^{-1}\in F_{2e}$.
\smallskip

The facts in the previous paragraphs imply that if $(A,1)\in K{\setminus}K_0$, then the matrix $AA^{\sigma}$ is diagonalisable over $F$ and $C$ is diagonalisable over $F_{2e}$. In particular, there exists a $P\in \GL(2,p^{2e})$ such that $P^{-1}CP=D'$ for $D'=\dia(\mu,\mu^{-1})$; multiplying by $\veps$ then gives $P^{-1}AA^{\sigma}P=D$ for $D=\dia(\lambda_1, \lambda_2)$. Here, either $\lambda_1,\lambda_2\in F_0$ with $\lambda_1\lambda_2\in N(F_0)$, or $\lambda_1,\lambda_2\in F{\setminus}F_0$ and $\lambda_1^{\sigma} = \lambda_2$. With $A$, $P$, $D$ and $D'$ as above, in $K$ we let $(B,1)=(P,0)^{-1}(A,1)(P,0) =(P^{-1}AP^{\sigma},1)$. Then,
$$(BB^{\sigma},0) = (P,0)^{-1}(A,1)(A,1)(P,0) = (P^{-1}AA^{\sigma}P,0) = (D,0)=(\veps D',0)$$
and it follows that $BB^{\sigma}=D= \veps D'$.
\smallskip

We now derive more details about the matrix $B=P^{-1}AP^{\sigma}$; recall that $P\in \GL(2,p^{2e})$. Let $u_1,u_2$ be linearly independent (column) eigenvectors of $C$ and $AA^{\sigma}$ for the eigenvalues $\mu,\mu^{-1}$ and $\lambda_1,\lambda_2$, respectively; we have $Cu_1=\mu u_1$, $Cu_2=\mu^{-1}u_2$, and  $AA^{\sigma}u_i= \lambda_iu_i$ for $i\in \{1,2\}$. Taking the $\sigma$-image of the last equation and then multiplying by $A$ from the left we obtain $AA^{\sigma}(Au_i^{\sigma})= \lambda_i^{\sigma}(Au_i^{\sigma})$ for $i=1,2$. This means that the column vectors $Au_i^{\sigma}$ are also eigenvectors of $AA^{\sigma}$ for the eigenvalues $\lambda_i^{\sigma}$, $i=1,2$. It follows that if $\lambda_i=\lambda_i^{\sigma}$ for $i=1,2$, then we must have $Au_i^{\sigma}=\veps_iu_i$, and if $\lambda_i=\lambda_{3-i}^{\sigma}$, then $Au_{i}^{\sigma}=\veps_{3-i}u_{3-i}$, in both cases for some $\veps_1,\veps_2\in F$. The last bit we need is the fact that for the matrix $P$ we may take $P=(u_1,u_2)$, i.e., the matrix formed by the columns $u_1,u_2$, with entries in $F_{2e}$. Now, for $i=1,2$, in the case $\lambda_i=\lambda_i^{\sigma}$ we have $AP^{\sigma}=(Au_1^{\sigma},Au_2^{\sigma})= (\veps_1u_1,\veps_2u_2) =P\dia(\veps_1,\veps_2)$, and in the case $\lambda_i=\lambda_{3-i}^{\sigma}$ a similar calculation gives $AP^{\sigma}=(Au_1^{\sigma},Au_2^{\sigma}) =(\veps_2u_2,\veps_1u_1)=P\off(\veps_1,\veps_2)$. This shows that our matrix $B=P^{-1}AP^{\sigma}$ is equal to $\dia(\veps_1,\veps_2)$ or to $\off(\veps_1,\veps_2)$ for suitable $\veps_1,\veps_2\in F$, depending on whether $\lambda_1^{\sigma}$ is equal to $\lambda_1$ or $\lambda_2$. In both cases, of course, $\veps_1\veps_2\in N(F)$.
\smallskip

Recalling our notation $[A,i]$ for the cosets $(A,i)L=\{(\delta A,i);\ \delta\in F^*\}$, the above calculations lead to the following result.
\smallskip

\bp\label{can1}
Let $G=M(p^{2f})$ for some odd prime $p$. Then, every element of the form $[A,1]\in G$ is conjugate in $\ovl{G}$ to $[B,1]$ with $B=\dia(\lambda,1)$ or $B=\off(\lambda,1)$ for some $\lambda\in N(F)$. If, in addition, $[AA^{\sigma},0]=[C,0]$ for some $C\in {\rm SL}(2,p^{2e})$ with $f/e$ odd, then $[B,1]=[P,0]^{-1}[A,1][P,0]$ for some $P\in {\rm GL}(2,p^{2e})$, and $\lambda\lambda^{\sigma}\in F_{2e}$ or $\lambda/\lambda^{\sigma}\in F_{2e}$, depending on whether $B$ is equal to $\dia(\lambda,1)$ or to $\off(\lambda,1)$.
\ep

\pr We have proven everything except for the last assertion. We have seen that if $[AA^{\sigma},0]=[C,0]$ for some $C\in {\rm SL}(2,p^{2e})$ with $f/e$ odd, then $BB^{\sigma}=P^{-1}(AA^{\sigma})P = \veps\dia(\mu,\mu^{-1})$ for some $\veps\in F^*$, $\mu\in F_{2e}$ and some $P\in {\rm GL}(2,p^{2e})$. If $B=\dia(\lambda,1)$, then we have  $\dia(\lambda\lambda^{\sigma}, 1)=BB^{\sigma}=\veps C=\veps\dia(\mu,\mu^{-1})$, which implies that $\veps=\mu$ and $\lambda\lambda^{\sigma}=\mu^2\in F_{2e}$. In the case when $B=\off(\lambda,1)$ we have $\off(\lambda,\lambda^{\sigma})=BB^{\sigma}=\veps C=\veps\dia(\mu,\mu^{-1})$, from which we obtain $\lambda/\lambda^{\sigma}=\mu^2\in F_{2e}$. \ebox

Let us have another look at conjugation in the group $\ovl{G}=J/L$. Observe that if $(P,i)\in J$, then $(P,i)^{-1} = ((P^{\sigma^i})^{-1},i)$. Conjugates of $(B,1)\in K$ by $(P,i)$ have the form  $(P,0)^{-1}(B,1)(P,0)= (P^{-1}BP^{\sigma},1)$ if $i=0$, and $(P,1)^{-1}(B,1)(P,1)=((P^{\sigma})^{-1}B^{\sigma}P,1)$ if $i=1$. It follows that two elements $(B,1)$ and $(B',1)$ of $K$ are conjugate in $J$ if and only if $B'=P^{-1}BP^{\sigma}$ or $B'=(P^{\sigma})^{-1}B^{\sigma}P$ for some $P\in \GL(2,F)$. Taking the $\sigma$-image in the second case and passing onto $G=K/L$ we have:
\smallskip

\bp\label{conj} Two elements $[B,1]$ and $[B',1]$ of $G$ are conjugate in $\ovl{G}$ if and only if $P^{-1}BP^{\sigma}=\veps B'$ or $P^{-1}BP^{\sigma}=\veps B'^{\sigma}$ for some $\veps \in F^*$ and some $P\in \GL(2,F)$. \hfill $\Box$
\ep

We will write the two conditions of Proposition \ref{conj} in the unified form $ P^{-1}BP^{\sigma}=
\veps B'^{(\sigma)} $, or, equivalently, $ BP^{\sigma}= \veps PB'^{(\sigma)} $ where $B'^{(\sigma)}$ is equal to $B'$ or $B'^{\sigma}$, depending on whether $i=0$ or $i=1$ when using the element $[P,i]$ for conjugation.
\smallskip

%------------------------------------------------
\section{Conjugacy of representatives of twisted elements}\label{sec:twi}
\smallskip

We continue with identification of elements of $\ovl{G}$ that conjugate a diagonal (or an off-diagonal) element from Proposition \ref{can1} to another such element. As a by-product we will be able to identify $\ovl{G}$-stabilisers of our representatives of twisted elements in $G$.
\smallskip

We begin with the case when $B=\dia(\lambda,1)$ and $B'=\dia(\lambda',1)$; by Proposition \ref{conj} it is sufficient to find the nonsingular matrices $P\in \GL(2,F)$ and $\veps\in F^*$ for which $ BP^{\sigma}= \veps PB'^{(\sigma)} $ in the sense of the notation introduced at the end of section \ref{sec:diag}. Throughout the computation we will use the symbols $\lambda^{(\sigma)}$ and $\lambda'^{(\sigma)}$ in an analogous way as explained for $B'^{(\sigma)}$. Assuming that $P$ has entries $\alpha,\beta,\gamma,\delta$, the above condition says that
$$ \begin{pmatrix} \lambda & 0 \\ 0 & 1\end{pmatrix}
\begin{pmatrix} \alpha^{\sigma} & \beta^{\sigma} \\ \gamma^{\sigma} & \delta^{\sigma} \end{pmatrix} = \veps \begin{pmatrix} \alpha & \beta \\ \gamma & \delta \end{pmatrix}
\begin{pmatrix} \lambda'^{(\sigma)} & 0 \\ 0 & 1\end{pmatrix} \ . $$
Evaluating the products we obtain:
$$\lambda \alpha^{\sigma}=\veps\lambda'^{(\sigma)}\alpha \ , \ \   \lambda\beta^{\sigma}=\veps\beta \ , \ \
\gamma^{\sigma}=\veps\lambda'^{(\sigma)}\gamma \ , \ \  \delta^{\sigma}=\veps\delta \ .$$
By straightforward manipulation this gives the following system of four equations:
$$ \alpha(\lambda\lambda^{\sigma} -\veps\veps^{\sigma}\lambda'\lambda'^{\sigma})=0 \ ,\ \
\beta(\lambda\lambda^{\sigma}-\veps\veps^{\sigma})=0\ ,\ \
\gamma(\veps\veps^{\sigma}\lambda'\lambda'^{\sigma}-1)=0 \ , \ \ \delta(\veps\veps^{\sigma}-1)=0\ .$$
From non-singularity of $P$ it follows that $\delta\ne 0$ or $\beta\ne 0$, that is, $\veps\veps^{\sigma}=1$ or $\veps\veps^{\sigma}=\lambda\lambda^{\sigma}$.
\smallskip

Consider first the case $\veps\veps^{\sigma}=1$, that is, $\veps^{q+1}=1$ for $q=p^f$. Since $\lambda,\lambda'\in N(F)$, we have $\lambda\lambda^{\sigma}\ne 1\ne \lambda'\lambda'^{\sigma}$. Our equations together with $\veps\veps^{\sigma}=1$ then imply that $\beta=\gamma=0$. Hence $\alpha,\delta\ne 0$, by non-singularity of $P$; in particular, $\lambda\lambda^{\sigma} = \lambda'\lambda'^{\sigma}$, or, equivalently, $(\lambda'/\lambda)^{q+1}=1$. We are interested in conjugation in the group $\ovl{G}=J/L$ and so we may assume that $\delta=1$, which reduces the relations below our matrix equation to $\veps=1$ and $\alpha^{q-1} = \lambda'^{(\sigma)}/\lambda$. Since $(\lambda'/\lambda)^{q+1}=1$, the equation $\eta^{q-1} = \lambda'/\lambda$ has $q-1$ solutions $\eta\in F^*$ (note that $|F^*|=q^2-1$). If $\lambda'^{(\sigma)}=\lambda'$, then all solutions of the equation $\alpha^{q-1} = \lambda'^{(\sigma)}/\lambda$ have the form $\alpha=\eta$, and if $\lambda'^{(\sigma)}= \lambda'^{\sigma}= \lambda'\lambda'^{q-1}$, then all solutions of this equation are $\alpha=\eta\lambda'$.
\smallskip

The second case to consider is $\veps\veps^{\sigma}=\lambda\lambda^{\sigma}\ (\ne 1)$, which implies that $\alpha=\delta=0$, and also $\lambda\lambda^{\sigma}\lambda'\lambda'^{\sigma}=1$ since $\gamma,\beta$ now must be non-zero. By the same token as above we may let $\gamma=1$ without loss of generality. Then, our equations for $\gamma$ and $\beta$ in this case reduce to $\veps\lambda'^{(\sigma)}=1$ and $\lambda\beta^{\sigma}=\veps\beta$, the latter now being equivalent to $\beta^{q-1} = 1/(\lambda\lambda'^{(\sigma)})$. Since now $(\lambda\lambda')^{q+1} = 1$, there are $q-1$ solutions $\zeta$ of the equation $\zeta^{q-1}=1/(\lambda\lambda')$ in $F^*$. If $\lambda'^{(\sigma)}= \lambda'$, then we have $\beta=\zeta$, and if $\lambda'^{(\sigma)}=\lambda'^{\sigma}= \lambda'\lambda'^{q-1}$, we have $\beta=\zeta/\lambda'$. Summing up, we arrive at the following:
\bigskip

\bp\label{p-dia} Let $B=\dia(\lambda,1)$ and $B'=\dia(\lambda',1)$ for $\lambda,\lambda'\in N(F)$. If an element $[P,i]\in \ovl{G}$ conjugates $[B,1]$ to $[B',1]$, then, without loss of generality, $P=\dia(\omega,1)$ or $P=\off(\omega,1)$ for suitable $\omega\in F^*$. Moreover:
\medskip

\noindent {\rm 1.} If $P=\dia(\omega,1)$, then $\lambda\lambda^{\sigma}= \lambda'\lambda'^{\sigma}$, and if this condition is satisfied, then $[B,1]$ conjugates to $[B',1]$ in $\ovl{G}$ exactly by the $q-1$ elements $[P,0]$ such that $\omega=\eta$ and the $q-1$ elements $[P,1]$ with $\omega=\eta\lambda'$, where $\eta\in F^*$ is one of the $q-1$ solutions of the equation $\eta^{q-1} = \lambda'/\lambda$.
\medskip

\noindent {\rm 2.} If $P=\off(\omega,1)$, then $\lambda\lambda^{\sigma} \lambda'\lambda'^{\sigma}=1$, and if this holds, then $[B,1]$ conjugates to $[B',1]$ in $\ovl{G}$ exactly by the $q-1$ elements $[P,0]$ with  $\omega=\zeta$ and the $q-1$ elements $[P,1]$ such that $\omega=\zeta/\lambda'$, where $\zeta\in F^*$ is one of the $q-1$ solutions of the equation $\zeta^{q-1}=1/(\lambda\lambda')$. \hfill $\Box$ \ep
\bigskip

We now repeat this process but now with matrices $B=\off(\lambda,1)$ and $B'=\off(\lambda',1)$. Conjugating by $[P,i]$ and assuming that $P$ has entries $\alpha,\beta,\gamma,\delta$, the unified form $ BP^{\sigma}= \veps PB'^{(\sigma)} $ of the condition of Proposition \ref{conj} now translates into the matrix equation
$$
\begin{pmatrix} 0 & \lambda \\ 1 & 0\end{pmatrix}
\begin{pmatrix} \alpha^{\sigma} & \beta^{\sigma} \\ \gamma^{\sigma} & \delta^{\sigma} \end{pmatrix} = \veps \begin{pmatrix} \alpha & \beta \\ \gamma & \delta \end{pmatrix}
\begin{pmatrix} 0 & \lambda'^{\sigma^i} \\ 1 & 0\end{pmatrix} \ . $$
It follows that
$$ \alpha^{\sigma}=\veps \delta\ , \ \ \ \beta^{\sigma}=\veps\lambda'^{\sigma^i}\gamma\ , \ \ \
\lambda\gamma^{\sigma} = \veps\beta\ , \ \ \ \lambda\delta^{\sigma} = \veps\lambda'^{\sigma^i}\alpha\ , $$
which, after some manipulation, yield the following two equations:
$$  \gamma(\lambda^{\sigma}-\veps\veps^{\sigma}\lambda'^{\sigma^i})=0 \ \ \ {\rm and} \ \ \ \ \
\delta^{\sigma}(\lambda - \veps\veps^{\sigma}\lambda'^{\sigma^i})=0\ .  $$
This all means that either (a) $\lambda = \veps\veps^{\sigma}\lambda'^{\sigma^i}$, and then $\beta=\gamma=0$ and we may assume $\delta=1$, or else (b) $\lambda^{\sigma} = \veps\veps^{\sigma}\lambda'^{\sigma^i}$, and then we have $\alpha=\delta=0$ and, without loss of generality, $\gamma=1$. Since $\veps\veps^{\sigma}, \lambda\lambda^{\sigma} \in F_0^*$, these conditions are equivalent to $\lambda/\lambda'\in F_0^*$ or $\lambda\lambda'\in F_0^*$, independently of the value of $i$, but in our analysis below it is still useful to refer to $i$.
\smallskip

In the case (a), when $\lambda/\lambda'^{\sigma^i} = \veps\veps^{\sigma}\in F_0^*$, for every $i\in \{0,1\}$ there are $q+1$ $(q+1)^{\rm th}$ roots $\eta_{(i)}$ of $\lambda/\lambda'^{\sigma^i}$ in $F^*$. From $\delta=1$ we have $a^{\sigma}=\veps$ and $\lambda=\veps\lambda'^{\sigma^i}\alpha$, that is, $\alpha^{q+1} = \lambda/\lambda'^{\sigma^i}$. This implies that $\alpha=\eta_{(i)}$ is one of the $(q+1)^{\rm th}$ roots of $\lambda/\lambda'^{\sigma^i}$, giving $q+1$ conjugation elements $[P,i]$ such that $P=\dia(\eta_{(i)},1)$. In the case (b), $\lambda^{\sigma}/\lambda'^{\sigma^i} = \veps\veps^{\sigma}\in F_0^*$ and since also $\lambda\lambda^{\sigma} \in F_0^*$, we have $\lambda\lambda'^{\sigma^i}\in F_0^*$. It follows that for every $i\in \{0,1\}$ there are $q+1$ $(q+1)^{\rm th}$ roots $\zeta_{(i)}$ of $\lambda\lambda'^{\sigma^i}$ in $F^*$. From $\gamma=1$ we obtain $\lambda=\veps\beta$ and $\beta^{\sigma}=\veps \lambda'^{\sigma^i}$, which means that $\beta^{q+1}=\lambda\lambda'^{\sigma^i}$. Consequently, $\beta=\zeta_{(i)}$ and we have in this second case $q+1$ conjugation elements $[P,i]$ such that $P=\off(\zeta_{(i)},1)$. Realising that the condition (a) for $i=0$ is equivalent to (b) for $i=1$ (and equivalent to $\lambda/\lambda'\in F_0^*$) and, similarly, the condition (a) for $i=1$ is equivalent to (b) for $i=0$ (and equivalent to $\lambda\lambda'\in F_0^*$), we conclude that:
\bigskip

\bp\label{p-off} Let $B=\off(\lambda,1)$ and $B'=\off(\lambda',1)$ for $\lambda,\lambda'\in N(F)$. Further, for $i\in \{0,1\}$, let $\eta_{(i)},\zeta_{(i)}\in F^*$ be any of the $q+1$ roots of the equation $\eta_{(i)}^{q+1}= \lambda/\lambda'^{\sigma^i}$ and $\zeta_{(i)}^{q+1}=\lambda\lambda'^{\sigma^i}$, respectively, Then, an element $[P,i]\in \ovl{G}$ conjugates $[B,1]$ to $[B',1]$ if and only if $\lambda/\lambda'^{\sigma^i}$ or $\lambda\lambda'^{\sigma^i}$ are elements of $F_0^*$. In an equivalent form, $[B,1]$ is conjugate to $[B',1]$ if and only if either
\medskip

\noindent {\rm 1.} $\lambda/\lambda'\in F_0^*$, in which case the conjugation is realised by exactly $q+1$ elements $[P,0]$ with $P=\dia(\eta_{(0)},1)$ and exactly $q+1$ elements $[P,1]$ with $P=\off(\zeta_{(1)},1)$, or
\medskip

\noindent {\rm 2.} $\lambda\lambda'\in F_0^*$, by the conjugation realised by exactly $q+1$ elements $[P,0]$ with $P=\dia(\eta_{(1)},1)$ and exactly $q+1$ elements $[P,1]$ with $P=\off(\zeta_{(0)},1)$. \hfill $\Box$ \ep

%------------------------------------------------
\section{Conjugacy classes of twisted elements}\label{sec:twist}
\smallskip

With the help of the calculations in the previous section we can now prove a useful result about identification of suitable representatives of conjugacy classes (in the group $\ovl{G}$) of elements of $G{\setminus}G_0$.
\smallskip

\bt\label{t-conj}
Let $\xi$ be a primitive element of $F$ and let $[A,1]$ be an element of $G$. Then, exactly one of the following two cases occur:
\begin{itemize}
\item[{\rm 1.}] There exists exactly one odd $i\in \{1,2,\ldots,(q-1)/2\}$ such that $[A,1]$ is conjugate in $\ovl{G}$ to $[B,1]$ with $B=\dia(\xi^i,1)$; the order of $[A,1]$ in $G$ is then $2(q-1)/\gcd\{q-1,i\}$.
\item[{\rm 2.}] There exists exactly one odd $i\in \{1,2,\ldots,(q+1)/2\}$ such that $[A,1]$ is conjugate in $\ovl{G}$ to $[B,1]$ with $B=\off(\xi^i,1)$, and the order of $[A,1]$ in $G$ is $2(q+1)/\gcd\{q+1,i\}$.
\end{itemize}
Furthermore, we have:
\begin{itemize}
\item[{\rm 3.}] The stabiliser of $[B,1]$ for $B=\dia(\lambda,1)$, $\lambda\in N(F)$, in the group $\ovl{G}$ is isomorphic to the cyclic group $Z_{2(q-1)}$ generated by (conjugation by) $[P,1]$ for $P=\dia(\mu\lambda,1)$ with a suitable $(q-1)^{\rm th}$ root of unity $\mu$, except when $\lambda$ is a $(q+1)^{\rm th}$ root of $-1$ and $q\equiv -1$ mod $4$, in which case the stabiliser is isomorphic to $Z_{2(q-1)}\cdot Z_2$.
\item[{\rm 4.}] The stabiliser of $[B,1]$ for $B=\off(\lambda,1)$, $\lambda\in N(F)$, in the group $\ovl{G}$ is isomorphic to the cyclic group $Z_{2(q+1)}$ generated by (conjugation by) $[P,1]$ for $P=\off(\mu\lambda,1)$, where $\mu$ is a suitable $(q+1)^{\rm th}$ root of unity, except when $\lambda$ is a $(q-1)^{\rm th}$ root of $-1$ and $q\equiv 1$ mod $4$, when the stabiliser is isomorphic to $Z_{2(q+1)}\cdot Z_2$.
\end{itemize}
\et

\pr
By Proposition \ref{can1}, each element $[A,1]\in G$ is conjugate in $G$ to $[B,1]$ with $B=\dia(\lambda,1)$ or $B=\off(\lambda,1)$ for some $\lambda\in N(F)$.
\smallskip

1. Suppose that $[A,1]$ is conjugate in $G$ to $[B,1]$, $B=\dia(\lambda,1)$, $\lambda\in N(F)$. Letting $B'=\dia(\lambda',1)$, Proposition \ref{p-dia} implies that the elements $[B,1]$ and $[B',1]$ are conjugate in $\ovl{G}$ if and only if $\lambda\lambda^{\sigma}= \lambda'\lambda'^{\sigma}$ or $\lambda\lambda^{\sigma} \lambda'\lambda'^{\sigma}=1$. If $\lambda=\xi^i$ and $\lambda'=\xi^{i'}$ for some odd $i$ and $i'$, the condition translates into $\xi^{(q+1)i}=\xi^{\pm(q+1)i'}$, which is equivalent to $(q+1)i\equiv \pm (q+1)i'$ mod $(q^2-1)$ and which simplifies to $i\equiv \pm i'$ mod $(q-1)$. This proves Part 1 except for the order assertion. But, in $G$, for $B=\dia(\xi^i,1)$, the order of $[B,1]$ is twice the order of $[B,1]^2=[BB^{\sigma},0]$, which is equal to the order of $BB^{\sigma}= \dia(\xi^{(q+1)i},1)$ in $PSL(2,q^2)$, which is known to be equal to $(q^2-1)/\gcd\{q^2-1,(q+1)i\} = (q-1)/\gcd\{q-1,i\}$.
\smallskip

2. Now, let $[A,1]$ be conjugate in $G$ to $[B,1]$ for $B=\off(\lambda,1)$ and $\lambda\in N(F)$. If   $B'=\off(\lambda',1)$, by Proposition \ref{p-off} we see that $[B,1]$ is conjugate to $[B',1]$ in $\ovl{G}$ if and only if $\lambda/\lambda'$ or $\lambda\lambda'$ are elements of $F_0^*$, or, equivalently, are powers of $\xi^{q+1}$. Letting $\lambda=\xi^i$ and $\lambda'=\xi^{i'}$ for odd $i,i'$, this condition translates into the congruence $i\equiv \pm i'$ mod $(q+1)$, proving Part 2 apart from the claim about the order. To fill in this last bit, the order of $[B,1]$ in $G$ for $B=\off(\xi^i,1)$ is two times the order of $[B,1]^2=[BB^{\sigma},0]$, which is equal to the order of $BB^{\sigma}=\dia(\xi^{i},\xi^{qi})$ in $PSL(2,q^2)$. The latter is equal to the order of $\dia(1,\xi^{(q-1)i})$ in $PSL(2,q^2)$, i.e., to $(q^2-1)/\gcd\{q^2-1,(q-1)i\} = (q+1)/\gcd\{q+1,i\}$.
\smallskip

3. The values $\lambda=\lambda'$ automatically satisfy the first condition in Proposition \ref{p-dia}, by which the stabiliser of $[B,1]$ for $B=\dia(\lambda,1)$ includes a subgroup $H$ formed by (conjugation by) the $(q-1)$ elements $[P,0]$ for $P=\dia(\eta,1)$ and the $q-1$ elements $[P,1]$ with $P=\dia(\eta\lambda,1)$, where $\eta$ ranges over the set of all $(q-1)^{\rm th}$ roots of unity in $F$ (and this set is simply equal to $F_0^*$). We claim that $H$ is (isomorphic to) a cyclic group of order $2(q-1)$. Indeed, let $\lambda=\xi^i$ for odd $i$, $1\le i\le (q-1)/2$. Then, letting $\eta = \xi^{(q+1)(1-i)/2}$, the group $H$ is generated by $[P,1]$, where $P= \dia(\eta\lambda)$. To see this, observe that $[P,1]^2=[PP^{\sigma},0]=[Q,0]$, where $Q=\dia(\eta^2\lambda\lambda^{\sigma},1)$. With $\lambda$ and $\eta$ as above we have $Q=\dia(\xi^{q+1},1)$, which means that the order of $[Q,0]$ is $q-1$ and so $[P,1]$ generates $H$, as claimed.
\smallskip

Observe that the values $\lambda=\lambda'$ may also satisfy our second condition $\lambda\lambda^{\sigma} \lambda'\lambda'^{\sigma}=1$ from Proposition \ref{p-dia}. This happens if and only if $\lambda$ is a $(q+1)^{\rm th}$ root of $-1$ in $F$ and $q\equiv -1$ mod $4$, which is what we now assume; note that now $\lambda^{\sigma}=-\lambda^{-1}$. Then, if $\zeta$ and $\zeta'$ are $(q-1)^{\rm th}$ roots of $(\lambda\lambda')^{-1}=\lambda^{-2}$ as in the calculations immediately preceding Proposition \ref{p-dia}, we have $\zeta'/\zeta \in F_0^*$, $\lambda^{q-1}=-\zeta^{q-1}$, and $[\off(\zeta'/\lambda,1),i] [\off(\zeta/\lambda,1),i]= [\dia((-1)^i\zeta'/\zeta,1),0]$ for $i\in\{0,1\}$. In particular, letting $Q=\off(\zeta/\lambda,1)$, the element $[Q,1]^2=[\dia(-1,1),0]$ lies in the cyclic group $H$; hence $[Q,1]$ has order $4$. Moreover, one may check that for every $\eta\in F_0^*$ we have $[Q,1]^{-1}[P,0][Q,1]=[P,0]^{-1}$ if $P=\dia(\eta,1)$ and $[Q,1]^{-1}[P,1][Q,1]=[P,1]^{-1}[Q,1]^2$ if $P=\dia(\eta\lambda,1)$.
\smallskip

Summing up, the stabiliser of an element $[B,1]$ with $B=\dia(\lambda,1)$ in $\ovl{G}$ for $\lambda\in N(F)$ is isomorphic to $H\cong Z_{2(q-1)}$ except for the case when $q\equiv -1$ mod $4$ and $\lambda$ is a $(q+1)^{\rm th}$ root of $-1$ and then the stabiliser of $[B,1]$ is isomorphic to an extension $H^*$ of $H$ by $[Q,1]$ with $[Q,1]^2\in H$, so that $H^*\cong Z_{2(q-1)}\cdot Z_2$, as claimed.
\smallskip

4. The values $\lambda=\lambda'$ automatically satisfy the condition (1) from Proposition \ref{p-off}. By the calculations preceding this Proposition, the stabiliser of $[B,1]$ for $B=\off(\lambda,1)$ contains a subgroup formed by the $(q+1)$ elements $[P,0]$ for $P=\dia(\eta,1)$ and the $q-1$ elements $[P,1]$ with $P=\off(\eta\lambda,1)$, where this time $\eta$ ranges over the set of all $(q+1)^{\rm th}$ roots of unity in $F$. The collection of these $2(q+1)$ elements is closed under multiplication and hence forms a group $H$ of order $2(q+1)$; we claim that $H$ is cyclic. To demonstrate this, let $\lambda=\xi^i$ for some odd $i$ such that $1\le i\le (q+1)/2$. Then, using $\eta = \xi^{(q-1)(1+i)/2}$, the element $[P,1]$ with $P= \off(\eta\lambda)$ is a generator of $H$. Indeed, a calculation shows that  $[P,1]^2=[PP^{\sigma},0]=[Q,0]$, where $Q=\dia(\eta^2\lambda^{1-q},1)$, and with the given $\lambda$ and our choice of $\eta$ we have $Q=\dia(\xi^{q-1},1)$. It follows that the order of $[Q,0]$ is $q+1$ and so $[P,1]$ generates $H$.
\smallskip

Here, the values $\lambda=\lambda'$ may also fulfil the condition (2) from Proposition \ref{p-off}, that is, $\lambda^2\in F_0^*$. Since we know that $\lambda\in N(F)$, it can be checked that $\lambda^2\in F_0^*$ is equivalent to $\lambda^{q-1}=-1$ and implies $q\equiv 1$ mod $4$. Again, the calculations made immediately before Proposition  \ref{p-off} imply the following: Letting $\zeta$ be any of the $q+1$ $(q+1)^{\rm st}$ roots of $\lambda^2$ in $F$, in the stabiliser of $[B,1]$ we have additional $q+1$ elements $[P,1]$ such that $P=\dia(\zeta/\lambda,1)$, as well as $q+1$ elements $[P,0]$ such that $P=\off(\zeta,1)$. Further calculations show that if $\zeta$ and $\zeta'$ are $(q+1)^{\rm th}$ roots of $\lambda^2$, $P=\off(\zeta,1)$ and $P'=\off(\zeta',1)$, then we have $[P',0][P,0]=[\dia(\zeta'/\zeta,1),0]\in H$ since it can be checked that $\zeta'/\zeta\in F_0^*$. Similarly, if $Q=\dia(\zeta/\lambda,1)$ and $Q'=\dia(\zeta'/\lambda,1)$, then $[Q',1][Q,1] = [Q'Q^{\sigma},0] = [\dia((\zeta'/\zeta)(\zeta/\lambda)^{q+1},1),0]\in H$. In particular letting $\zeta'=\zeta$ we have $(\zeta/\lambda)^{q+1}=\lambda^2/(\lambda^{q-1}\lambda^2)=-1$ and so $[Q,1]^2=[\dia(-1,1),0]\in H$. We thus have a situation analogous to the Case 3 of this proof: $[Q,1]$ has order $4$, conjugation of an element $[P,0]\in H$ with $P=\dia(\eta,1)$ by $[Q,1]$ inverts $[P,0]$, while conjugation of a $[P,1]\in H$ such that  $P=\off(\eta\lambda,1)$ by $[Q,1]$ gives $[P,1]^{-1}[Q,1]^2$.
\smallskip

In summary, the stabiliser of an element $[B,1]$ for $B=\off(\lambda,1)$ is isomorphic to $H\cong Z_{2(q+1)}$ if $\lambda$ is not a $(q-1)^{\rm th}$ root of $-1$, and isomorphic to an extension $H^*$ of $H$ by $[Q,1]$ (where  $[Q,1]^2\in H$), implying that $H^*\cong Z_{2(q-1)}\cdot Z_2$. This completes the proof.
\ebox

Let us remark that the exceptional cases in the items 3 and 4 above correspond precisely to elements $[B,1]$ of order $4$ in $G$. By inspecting possible orders of $[B,1]$ we also have:

\bc\label{c-conj}
Every element of $G{\setminus} G_0$ has order divisible by $4$. \hfill $\Box$
\ec
\smallskip

%------------------------------------------------
\section{Non-singular pairs and twisted subgroups}\label{sec:twpairs}
\smallskip

Our aim in this and the following two sections is to determine representatives of selected conjugacy classes $\{(x,y)^g;\ g\in \ovl{G}=G\langle\sigma\rangle\}$ of elements $x,y\in G$ satisfying $(xy)^2=1$, and make important conclusions about subgroups the corresponding pairs $(x,y)$ generate. We note that the action of $\sigma$ need not be considered separately, because $[I,1][A,i][I,1]= [A^{\sigma},i]$ for $i\in \{0,1\}$, which means that the action of $\sigma$ is equivalent to conjugation by the element $[I,1]\in \ovl{G}$.
\smallskip

Since we want $xy$ to have order $2$, both $x$ and $y$ as above must lie in $G{\setminus}G_0$ because, by Corollary \ref{c-conj}, there are no involutions in $G{\setminus}G_0$. By Theorem \ref{t-conj} we may assume that $y=[B,1]$ for $B=\dia(\lambda,1)$ or $B=\off(\lambda,1)$ for a suitable $\lambda\in N(F)$. Letting $x=[A,1]$, the pair $x,y$ may in general generate a proper subgroup of $G=M(p^{2f})$; such cases will still be of interest for our intended classification of orientably-regular maps as long as the subgroup $\langle x,y\rangle$ is twisted, that is, isomorphic to $M(p^{2e})$ for a suitable divisor $e$ of $f$.
\smallskip

We now identify conditions on $A$ implied by the requirement that $([A,1][B,1])^2$ be the identity in $G$ and begin with the case when $B=\dia(\lambda,1)$. Let $A\in \GL(2,F)$ have rows $(a,b)$ and $(c,d)$, with determinant $ad-bc\in N(F)$. Then, $[A,1][B,1] = [AB^{\sigma},0]$, where
\[ AB^{\sigma}=\begin{pmatrix} a\lambda^{\sigma} & b \\ c\lambda^{\sigma} & d\end{pmatrix} \ .\] Since $AB^{\sigma}$ lies in $\PSL(2,F)$, it has order $2$ if and only if its trace $a\lambda^{\sigma}+d$ is equal to zero. If one of $a,d$ was equal to zero, both would have to be zero and then $[A,1]$ and $[B,1]$ would clearly not generate a twisted subgroup of $G$. Therefore both $a,d$ are non-zero and we may assume without loss of generality that $a=-1$ and $d=\lambda^{\sigma}$. We will thus consider only elements $[A,1]\in G$ of the form
\begin{equation}\label{A1} A=\begin{pmatrix} -1 & b \\ c & \lambda^{\sigma}\end{pmatrix}\ , \ \ u=bc + \lambda^{\sigma}\in N(F)\ .\end{equation}

Next, consider the case when $B=\off(\lambda,1)$. For a matrix $A\in \GL(2,F)$ with rows $(a,b)$ and $(c,d)$ such that $ad-bc\in N(F)$ we now have $[A,1][B,1] = [AB^{\sigma},0]$, where \[ AB^{\sigma}= \begin{pmatrix} b & a\lambda^{\sigma} \\ d & c\lambda^{\sigma}\end{pmatrix}\ .\] Again, $AB^{\sigma}\in \PSL(2,F)$ has order $2$ if and only if its trace $b+c\lambda^{\sigma}$ is equal to zero. If one of $b,c$ was equal to zero, we would have $b=c=0$, but then $[A,1]$ and $[B,1]$ would again not generate a twisted subgroup of $G$. Therefore both $b$ and $c$ are non-zero and we may assume that $c=-1$ and $b=\lambda^{\sigma}$. It follows that, without loss of generality, we only need to consider elements $[A,1]\in G$ such that \begin{equation}\label{A2} A=\begin{pmatrix} a & \lambda^{\sigma} \\ -1 & d \end{pmatrix} \ , \ \ u=ad + \lambda^{\sigma}\in N(F)\ .\end{equation}

With $A$ and $B$ as above we can still identify obvious instances when $[A,1]$ and $[B,1]$ do not generate a twisted subgroup of $G$. This is certainly the case if
\vskip -3mm

\begin{itemize}
\item[{\rm (i)}] both $[A,1]$ and $[B,1]$ have order $4$, as then the two elements generate a solvable group, cf. \cite{CoMo}, or
\item[{\rm (ii)}] $B=\dia(\lambda,1)$ and $A$ is an upper- or a lower-triangular matrix, as then $[A,1]$ and $[B,1])$ generate a triangular subgroup of $G$, or else
\item[{\rm (iii)}] $B=\off(\lambda,1)$ and $A$ is an off-diagonal matrix, as then $[A,1]$ and $[B,1]$ clearly do not generate a twisted subgroup of $G$.
\end{itemize}
\vskip -2mm

\noindent For $B=\dia(\lambda,1)$ and $A$ given by (\ref{A1}) and for $B=\off(\lambda,1)$ and $A$ given by (\ref{A2}), an ordered pair $([A,1],[B,1])$ not satisfying any of (i), (ii) and (iii) will be called {\em non-singular}.
\smallskip

We are now in position to classify the subgroups of $G=M(p^{2f})$ generated by non-singular pairs. To do so we will again use knowledge of the situation in the subgroup $G_0\simeq {\rm PSL}(2,p^{2f})$ of $G$. Recall that a subgroup $H$ of $G=M(p^{2f})$ was said to stabilise a point if there exists an element in $F\cup \{\infty\}$ fixed by all linear fractional mappings corresponding to elements of $H$; also, $G_0$ denotes the (unique) copy of ${\rm PSL}(2,p^{2f})$ in $G$.
\smallskip

\bp\label{non-sin}
Let $H$ be a subgroup of $G$ generated by a non-singular pair $([A,1],[B,1])$. Then $H$ is isomorphic to $M(p^{2e})$ for some positive divisor $e$ of $f$ such that $f/e$ is odd.
\ep

\pr
Let $H_0=H\cap G_0$. The classification of \cite{LS} tells us that subgroups of $G_0$ fall into the following categories: point stabilisers, dihedral groups, $A_4$, $S_4$, $A_5$, and ${\rm PSL}(2,p^{g})$ for divisors $g$ of $2f$ together with ${\rm PGL}(2,p^{g/2})$ for even divisors $g$ of $2f$. For our subgroup $H_0$ we subsequently rule out all but the penultimate case.
\smallskip

We first show that $H_0$ neither stabilises a point nor is dihedral. The first coordinate of every element of the group $H_0$ can be written as a product of the form $X_1Y_1^{\sigma}\ldots X_mY_m^{\sigma}$, where $X_i,Y_i\in \{A,B\}$, $1\le i\le m$. It follows that $H_0$ is generated by the three elements $[BB^{\sigma},0]$, $[BA^{\sigma},0]$ and $[AB^{\sigma},0]$; note that $AA^{\sigma}=AB^{\sigma}(BB^{\sigma})^{-1} BA^{\sigma}$. A straightforward calculation (details of which we leave to the reader) shows that for our generators $[A,1]$ and $[B,1]$ of $H$, there is no point in $F\cup \{\infty\}$ fixed by all three of the linear fractional mappings corresponding to $BB^{\sigma}$, $BA^{\sigma}$ and $AB^{\sigma}$. It follows that the subgroup $H_0=\langle [BB^{\sigma},0],[BA^{\sigma},0],[AB^{\sigma},0]\rangle$ cannot stabilise a point.
\smallskip

Suppose that $H_0$ is dihedral. An easy calculation shows that we always have $AB^{\sigma}\ne BA^{\sigma}$. It follows that one of $AB^{\sigma}$, $BA^{\sigma}$ is an involution outside the cyclic part of $H_0$ and hence inverts one of $AA^{\sigma}$ and $BB^{\sigma}$. If conjugation by $AB^{\sigma}$ is considered, the conditions $(B^{\sigma})^{-1}A^{-1}XX^{\sigma}AB^{\sigma} =(X^{\sigma})^{-1}X^{-1}$ for $X\in \{A,B\}$ both reduce to $(B^{-1}A)^2=I$ in $\PSL(2,p^{2f})$. A similar calculation shows that the condition of $XX^{\sigma}$ being inverted by conjugation by $BA^{\sigma}$ for $X\in \{A,B\}$ both reduce to $(B^{-1}A)^2=I$ again. The condition $(B^{-1}A)^2=I$ in $\PSL(2,p^{2f})$ is equivalent to $(B^{-1}A)^2=\veps I$ in ${\rm SL}(2,p^{2f})$ for some $\veps\in {\rm GF}(p^{2f})$. If $B=\dia(\lambda,1)$ and $b=c=0$ in $A$, then $\lambda\lambda^{\sigma}=-1$ which, by part 3 of Theorem \ref{t-conj}, means that the orders of both $[A,1]$ and $[B,1]$ would be equal to $4$. If one of $b,c\ne 0$, then $\veps=-1$ and $\lambda\lambda^{\sigma}=1$, a contradiction as $\lambda\in N(F)$. An entirely similar conclusion is obtained if $B=\off(\lambda,1)$: if $a=d=0$ in $A$, then $\lambda^{\sigma}/\lambda=-1$ and $[A,1]$, $[B,1]$ have order $4$ by part 4 of Theorem \ref{t-conj}, and if one of $a,d\ne 0$, then $\veps=-1$ and $\lambda^{\sigma}=\lambda$, a contradiction again. It follows that for all non-singular pairs $([A,1],[B,1])$ we have $(B^{-1}A)^2\ne I$ in $\PSL(2,p^{2f})$, and so $H_0$ cannot be dihedral.
\smallskip

By Corollary \ref{c-conj}, the order of every element in $G{\setminus} G_0$ is divisible by $4$. Since the order of one of $[A,1]$ and $[B,1]$ is assumed to be greater than $4$, one of the elements $[AA^{\sigma},0]$, $[BB^{\sigma},0]$ of $H_0$ has to have even order greater than $2$. But as $A_4$ and $A_5$ do not contain elements of even order greater than $2$, $H_0$ cannot be isomorphic to $A_4$ or $A_5$. The case $H_0\cong S_4$ is excluded as follows. If there was an orientably-regular map with automorphism group $H$ such that $H_0\cong S_4$, then the orders of the elements $[A,1]$, $[B,1]$ would have to be in the set $\{4,8\}$ since the largest even order of an element in $S_4$ is $4$. According to the list of \cite{CoList}, however, there is no orientably-regular map of type $\{4,8\}$ or $\{8,8\}$ with (orientation-preserving) automorphism group of order $2|S_4|=48$.
\smallskip

We may now apply Proposition \ref{subgps-class} to the subgroup $H_0$ to conclude that $H$ is conjugate in $\ovl{G}$ to a subgroup isomorphic to $M(p^{2e})$ for some positive divisor $e$ of $f$ such that $f/e$ is odd, completing the proof.
\ebox

It follows that a pair $([A,1],[B,1])$ of elements of $G$ as above generates a twisted subgroup of $G$ if and only if the pair is non-singular.
\smallskip

%------------------------------------------------
\section{Orbits of non-singular pairs: The diagonal case} \label{sec:rep-dia}
\smallskip

We will identify representatives of $\ovl{G}$-orbits of non-singular pairs $([A,1],[B,1])$, dealing with $B=\dia(\lambda,1)$ and $A$ as in (\ref{A1}) here and deferring the case $B=\off(\lambda,1)$ to the next section.
\smallskip

Instead of working with matrices, the form of $A$ in (\ref{A1}) suggests to look at the corresponding quadruples $(\lambda,b,c,u)$, also called {\em non-singular}, under the induced action of the stabiliser of $[B,1]$ in $\ovl{G}$. We recall that the values of $\lambda$ and identification of the stabiliser are in items 1 and 3 of Theorem \ref{t-conj}. To simplify the notation in what follows, for any $\omega\in F$ we will use the symbol $\sqrt[r]{\omega}$ to denote the set of all $r^{\rm th}$ roots of $\omega$ in $F=\GF(q^2)$, $q=p^f$. The analysis in the third part of the proof of Theorem \ref{t-conj} tells us that the stabiliser of $[B,1]$ in $\ovl{G}$ consists exactly of the following elements of $\ovl{G}$:
\smallskip

\begin{tabular}{ll}
$[P_1(\eta),0]$, & where $P_1=\dia(\eta,1)$ and $\eta\in F_0^*$; \\
$[P_2(\eta),1]$, & where $P_2=\dia(\eta\lambda,1)$ and $\eta\in F_0^*$; \\
$[P_3(\zeta),0]$, & where $P_3=\off(\zeta,1)$ if $\lambda\in \sqrt[q{+}1]{-1}$ and $\zeta\in \sqrt[q{-}1]{\lambda^{-2}}$; \\
$[P_4(\zeta),1]$, & where $P_4=\off(\zeta/\lambda,1)$ if $\lambda\in \sqrt[q{+}1]{-1}$ and $\zeta\in \sqrt[q{-}1]{\lambda^{-2}}$.
\end{tabular}
\smallskip

\noi To find the corresponding orbit of $[A,1]$ we first evaluate the products $[P_j(\eta),0]^{-1}[A,1][P_j(\eta),0]$ for $A$ as in (\ref{A1}) and $j\in \{1,2,3,4\}$:
\smallskip

\begin{tabular}{ll}
$[P_1(\eta),0]^{-1}[A,1][P_1(\eta),0] = [C_1,1]$, &  where \ \  $C_1=\begin{pmatrix} -1 & b\eta^{-1} \\ c\eta &  \lambda^{\sigma}\end{pmatrix};$ \\
$[P_2(\eta),1]^{-1}[A,1][P_2(\eta),1] = [C_2,1]$, &  where \ \ $C_2=\begin{pmatrix} -1, & b^{\sigma}(\eta\lambda)^{-1} \\ c^{\sigma}\eta\lambda^{\sigma} &  \lambda^{\sigma}\end{pmatrix};$ \\
$[P_3(\zeta),0]^{-1}[A,1][P_3(\zeta),0] = [C_3,1]$, &  where \ \ $C_3=\begin{pmatrix} -1 & c\zeta/\lambda \\ b\lambda/\zeta & \lambda^{\sigma}\end{pmatrix};$ \ \ and \\
$[P_4(\zeta),1]^{-1}[A,1][P_4(\zeta),1] = [C_4,1]$, &  where \ \  $C_4=\begin{pmatrix} -1 &  -c^{\sigma}\zeta/\lambda^2 \\  b^{\sigma}/\zeta & \lambda^{\sigma}\end{pmatrix}.$
\end{tabular}

Let $\lambda=\xi^i$ for a fixed primitive element $\xi\in F$ and some odd $i$ such that $1\le i\le (q-1)/2$; note that here $\lambda\in \sqrt[q{+}1]{-1}$ if and only if $i=(q-1)/2$. It follows that we have either $(q-1)/4$ such odd values of $i$ if $q\equiv 1$ mod $4$ and all are smaller than $(q-1)/2$, or else $(q-3)/4$ such odd $i<(q-1)/2$ together with $i=(q-1)/2$ if $q\equiv -1$ mod $4$.
\smallskip

In each of the above cases, if $i<(q-1)/2$, the stabiliser $H$ of $[B,1]$ for $B=\dia(\lambda,1)$ has order $2(q-1)$. By the above calculations leading to the matrices $C_1$ and $C_2$, the orbit $O$ of the induced action of $H$ on a quadruple $(\lambda,b,c,u)$ satisfying $bc+\lambda^{\sigma}=u\in N(F)$ is formed by the quadruples $(\lambda,b\eta^{-1},c\eta,u)$ for $\eta\in F_0^*$ and $(\lambda, b^{\sigma}(\eta\lambda)^{-1}, c^{\sigma}\eta\lambda^{\sigma}, u^{\sigma}\lambda^{\sigma} /\lambda)$ for $\eta\in F_0^*$. It is easy to check that all these $2(q-1)$ quadruples are mutually distinct and hence each such orbit $O$ has size $2(q-1)$. If $u=\lambda^{\sigma}$, then one of $b,c$ would have to be equal to zero, contrary to non-singularity. For each of our $\lfloor(q-1)/4\rfloor$ choices of $i$ we therefore have a total of $(q^2-3)/2$ choices for $u\in N(F)$, $u\ne \lambda^{\sigma}$; for each of these we can freely choose a non-zero $b\in F$ in $q^2-1$ ways and then $c$ is uniquely determined. This gives a total of $\lfloor(q-1)/4\rfloor(q^2-1)(q^2-3)/2$ quadruples $(\lambda,b,c,u)$ and the number of $H$-orbits on these is given by dividing by $|H|=2(q-1)$. Denoting by $n_1$ the number of $\ovl{G}$-equivalence classes of non-singular quadruples $(\lambda,b,c,u)$ for $i<(q-1)/2$ for diagonal $B$, we obtain
\begin{equation}\label{bc1} n_1=(q+1)\left\lfloor\frac{q-1}{4}\right\rfloor\frac{q^2-3}{4} \ .\end{equation}

The enumeration is trickier for $i=(q-1)/2$, that is, for $\lambda=\xi^{(q-1)/2}$ and $q\equiv -1$ mod $4$. This is the case when the stabiliser $H^*$ of $[B,1]$ for $B=\dia(\lambda,1)$ has order $4(q-1)$ and the orbit $O^*$ of the induced action of $H^*$ on a quadruple $(\lambda,b,c,u)$ satisfying $bc+\lambda^{\sigma}=u\in N(F)$ includes: (i) the $2(q-1)$ quadruples in $O$ listed in the previous paragraph, and (ii) by the form of the matrices $C_3$ and $C_4$ also the $(q-1)$ quadruples $(\lambda,b',c',u)$ for $b'=c\zeta/\lambda$, $c'=b\lambda/\zeta$ as well as the $(q-1)$ quadruples $(\lambda,b'',c'',u/\lambda^2)$ for $b'' = -c^{\sigma}\zeta/\lambda^2$, $c''= b^{\sigma}/\zeta$, where $\zeta$ can be equal to any of the $q-1$ values of $\sqrt[q{-}1] {\lambda^{-2}}$. Clearly, either $|O^*|=4(q-1)$ or $O^*=O$, and we will investigate the second possibility.
\smallskip

We have $O^*=O$ if and only if $b$ is equal to one of $b'$ or $b''$. For our $\lambda =\xi^{(q-1)/2}$ we find that $\sqrt[q{-}1] {\lambda^{-2}} = \{\xi^{-1}\omega;\ \omega\in F_0^*\}$. It can be checked that $b$ cannot be equal to any of the values of $b''$ but it may happen that $b=b'= c\zeta/\lambda=c\xi^{-1}\xi^{(1-q)/2}\omega'$ for some $\omega'\in F_0^*$.  Taking into account the fact that we are only interested in non-zero $b,c$ it follows that $O^*=O$ if and only if $c/b=\xi^{(q+1)/2}\omega$ for some $\omega\in F_0^*$. But this is equivalent to $(c/b)^{q-1}=-1$, which can be rewritten in the form $bc^{\sigma}+b^{\sigma} c=0$. Since $bc^{\sigma}+b^{\sigma} c$ is the trace of $AA^{\sigma}$ for $A$ as in (\ref{A1}), the fact that the trace is zero is equivalent to $[A,1]^2$ having order $2$ in $G$. This all means that $O^*=O$ if and only if the element $[A,1]$ has order $4$ in $G$. But for $\lambda =\xi^{(q-1)/2}$ the order of $[B,1]$ in $G$ is then equal to $4$ as well, a contradiction with non-singularity. It follows that $|O^*|=4(q-1)$.
\smallskip

Observe that if $c/b=\xi^{(q+1)/2}\omega$ for some $\omega\in F_0^*$, then the determinant equation gives  $u=bc+\lambda^{\sigma}=b^2\xi^{(q+1)/2}\omega + \lambda^{\sigma}$, and since $q\equiv -1$ mod $4$ this means that $u-\lambda^{\sigma} \in S(F)$. Conversely, if this holds, then for every $\omega\in F_0^*$ the equation $b^2\xi^{(q+1)/2}\omega=u - \lambda^{\sigma}$ has two values of $b$ as a solution. For a fixed $v\in N(F)$ let $N(v)$ denote the set of all  $u\in N(F){\setminus}\{v\}$ for which $u-v\in S(F)$. It is easy to see that the value of $|N(v)|$ does not depend on the choice of $v\in N(F)$. Let $n_F$ denote this common value; it was shown in \cite{Kel} that $n_F=(q^2-1)/4$.
\smallskip

It follows that for every $u\in N(\lambda^{\sigma})$ we have $2(q-1)$ values of $b$ satisfying $b^2\xi^{(q+1)/2}\omega=u - \lambda^{\sigma}$ for some $\omega\in F_0^*$, and these yield the singular pairs identified in the previous paragraph. Thus, for any of the $n_F$ values of $u\in N(\lambda^{\sigma})$ we have $(q^2-1)-2(q-1)=(q-1)^2$ values of $b\in F^*$ that furnish a non-singular pair $[A,1],[B,1]$ for $G$, giving a total of $n_F(q-1)^2$ quadruples of the form $(\xi^{(q-1)/2},b,c,u)$. Any of the remaining $(q^2-3)/2 - n_F$ non-square values of $u$ then simply give $((q^2-3)/2 - n_F)(q^2-1)$ non-singular quadruples $(\xi^{(q-1)/2},b,c,u)$ that we have to consider. As established earlier, the induced action of the group $H^*$ on these quadruples is semi-regular. Hence, letting $n_2$ be the number of $\ovl{G}$-equivalence classes of non-singular quadruples $(\lambda,b,c,u)$ for $\lambda=\xi^{(q-1)/2}$ and $q\equiv -1$ mod $4$ and for diagonal $B$, we have
\begin{equation}\label{bc2} n_2=  \frac{n_F(q-1)^2 + ((q^2-3)/2 - n_F)(q^2-1)}{4(q-1)} = \frac{1}{8}\left((q+1)(q^2-3) - 4n_F\right)\ .\end{equation}

%------------------------------------------------
\section{Orbits of non-singular pairs: The off-diagonal case} \label{sec:rep-off}
\smallskip

We continue looking for representatives of $\ovl{G}$-orbits of non-singular pairs $([A,1],[B,1])$ as in the previous section but this time for $B=\off(\lambda,1)$ and $A$ as in (\ref{A2}). We again reduce the task to finding representatives of the corresponding quadruples $(\lambda,a,d,u)$, called {\em non-singular} as well, under the induced action of the stabiliser of $[B,1]$ in $\ovl{G}$, with $\lambda$ and the stabiliser being determined by items 2 and 4 of Theorem \ref{t-conj}. The orbit of $[A,1]$ given by (\ref{A2}) under conjugation by elements stabilising $[B,1]$ is found by calculating the conjugates by the following elements of $\ovl{G}$ identified in the fourth part of the proof of Theorem \ref{t-conj}:
\smallskip

\begin{tabular}{ll}
$[P_1(\eta),0]$, & where $P_1=\dia(\eta,1)$ and $\eta\in \sqrt[{q{+}1}]{1}$; \\
$[P_2(\eta),1]$, & where $P_2=\off(\eta\lambda,1)$ and $\eta\in \sqrt[{q{+}1}]{1}$; \\
$[P_3(\zeta),0]$, & where $P_3=\off(\zeta,1)$ if $\lambda\in \sqrt[q{-}1]{-1}$ and $\zeta\in \sqrt[{q{+}1}]{\lambda^2}$; \\
$[P_4(\zeta),1]$, & where $P_4=\dia(\zeta/\lambda,1)$ if $\lambda\in \sqrt[q{-}1]{-1}$ and $\zeta\in \sqrt[{q{+}1}]{\lambda^2}$.
\end{tabular}
\smallskip

\noi Evaluation of the products $[P_j(\eta),0]^{-1}[A,1][P_j(\eta),0]$ for $A$ as in (\ref{A2}) and $j\in \{1,2,3,4\}$ gives:
\smallskip

\begin{tabular}{ll}
$[P_1(\eta),0]^{-1}[A,1][P_1(\eta),0] = [D_1,1]$, & where \ \  $D_1=\begin{pmatrix} a\eta^{-1} & \lambda^{\sigma} \\ -1 & d\eta \end{pmatrix}$; \\
$[P_2(\eta),1]^{-1}[A,1][P_2(\eta),1] = [D_2,1]$, & where \ \ $D_2=\begin{pmatrix} -\eta^{-1}d^{\sigma}\lambda^{\sigma}/\lambda & \lambda^{\sigma} \\ -1 & -\eta a^{\sigma} \end{pmatrix}$; \\
$[P_3(\zeta),0]^{-1}[A,1][P_3(\zeta),0] = [D_3,1]$, & where \ \ $D_3=\begin{pmatrix} d\zeta/\lambda &  \lambda^{\sigma} \\ -1 & a\lambda/\zeta \end{pmatrix}$; \ and \\
$[P_4(\zeta),1]^{-1}[A,1][P_4(\zeta),1] = [D_4,1]$, & where \ \ $D_4=\begin{pmatrix}-a^{\sigma}\zeta/\lambda & \lambda^{\sigma} \\ -1 & d^{\sigma}\lambda/\zeta \end{pmatrix}$\ .
\end{tabular}

Let $\lambda=\xi^i$ for a fixed primitive element $\xi\in F$ and some odd $i$ such that $1\le i\le (q+1)/2$; here we have $\lambda\in \sqrt[q{-}1]{-1}$ if and only if $i=(q+1)/2$. This gives either $(q+1)/4$ such odd values of $i$ if $q\equiv -1$ mod $4$, all smaller than $(q+1)/2$, or $(q-1)/4$ such odd $i<(q+1)/2$ together with $i=(q+1)/2$ if $q\equiv 1$ mod $4$.
\smallskip

In each of these instances, if $i<(q+1)/2$, the stabiliser $H$ of $[B,1]$ for $B=\off(\lambda,1)$ has order $2(q+1)$. By the above calculations the orbit $O$ of the induced action of $H$ on a quadruple $(\lambda,a,d,u)$ satisfying $ad+\lambda^{\sigma}=u\in N(F)$ consists of the quadruples $(\lambda,a\eta^{-1},d\eta,u)$ and $(\lambda, -d^{\sigma}\lambda^{\sigma} (\eta\lambda)^{-1}, -a^{\sigma}\eta, u^{\sigma}\lambda^{\sigma} /\lambda)$ for $\eta\in F_0^*$. By non-singularity, at least one of $a,d$ is non-zero. This implies, by an easy computation, that all the above $2(q+1)$ quadruples are mutually distinct and hence each such orbit $O$ on our quadruples under the action of $H$ has size $2(q+1)$. If $u=\lambda^{\sigma}$, then exactly one of $a,d$ is equal to zero while the other one can assume any of the $q^2-1$ values in $F^*$. For each of the $\lfloor(q+1)/4\rfloor$ choices of $i$ considered here, we have:
\smallskip

($\alpha$) \  $(q^2-3)/2$ choices of $u\in N(F){\setminus}\{\lambda^{\sigma}\}$, each admitting a choice of $a\in F^*$ in $q^2-1$ ways, uniquely determining $d$ and making $\lfloor(q+1)/4\rfloor(q^2-1)(q^2-3)/2$ quadruples; and
\smallskip

($\beta$) \ for $u=\lambda^{\sigma}$ we have $2(q^2-1)$ choices of non-zero $a$ (for $d=0$) and non-zero $d$ (for $a=0$), giving $\lfloor(q+1)/4\rfloor\cdot 2(q^2-1)$ quadruples.
\smallskip

\noindent Since the number of $H$-orbits on quadruples in ($\alpha$) and ($\beta$) is given by dividing by $|H|=2(q+1)$, we arrive at a total of $n_3$ $\ovl{G}$-equivalence classes of non-singular quadruples $(\lambda,a,d,u)$ for $i<(q+1)/2$ and off-diagonal $B$, where
\begin{equation}\label{ad1} n_3=(q-1)\left\lfloor\frac{q+1}{4}\right\rfloor\frac{q^2+1}{4} \ . \end{equation}

The enumeration is again more tricky for $i=(q+1)/2$, that is, for $\lambda=\xi^{(q+1)/2}$ and $q\equiv 1$ mod $4$. The stabiliser $H^*$ of $[B,1]$ for $B=\off(\lambda,1)$ has in this case order $4(q+1)$ and the orbit $O^*$ of the induced action of $H^*$ on a quadruple $(\lambda,a,d,u)$ with values fulfilling the equation $ad+\lambda^{\sigma} =u\in N(F)$ contains the $2(q+1)$ quadruples in $O$ listed in the previous paragraph together with the $2(q+1)$ distinct quadruples $(\lambda,d\zeta/\lambda,a\lambda/\zeta,u)$ and $(\lambda,-a^{\sigma}\zeta/\lambda, d^{\sigma}\lambda/\zeta,-u^{\sigma})$, where $\zeta$ can be any of the $q+1$ values of $\sqrt[q{+}1] {\lambda^{2}}$. As before, either $|O^*|=4(q+1)$ or $O^*=O$, and we will study the second eventuality in detail.
\smallskip

By inspection of the orbits $O$ and $O^*$ in this case, the orbit $O^*$ coincides with $O$ if and only if either
\smallskip

($\gamma$) \ $a\eta^{-1}=d\zeta/\lambda$, or
\smallskip

($\delta$) \ $a\eta^{-1}=-a^{\sigma}\zeta/\lambda$ and $d\eta=d^{\sigma}\lambda/\zeta$,
\smallskip

\noindent all for some $\eta\in \sqrt[q{+}1] {1}$. For our calculations we may without loss of generality assume that $\zeta=\xi$, as this is one of the $(q+1)^{\rm th}$ roots of $\lambda^2=\xi^{q+1}$. Note that we now have $\lambda^{\sigma} = -\lambda$, $\eta^{\sigma}=\eta^{-1}$, $\zeta^{\sigma}=\lambda^2/\zeta$, and therefore also $(\zeta\eta/\lambda)^{\sigma} = -(\zeta\eta/\lambda)^{-1}$.
\smallskip

The case ($\gamma$): From $a=d(\zeta\eta/\lambda)$ and its $\sigma$-image $a^{\sigma}=d^{\sigma}(\zeta\eta/\lambda)^{\sigma} = -d^{\sigma}(\zeta\eta/\lambda)^{-1}$ we obtain $aa^{\sigma}+dd^{\sigma}=0$. But $aa^{\sigma}+dd^{\sigma}$ is the trace of $AA^{\sigma}$, which means that the order of both $[A,1]$ and $[B,1]$ is $4$. It follows that ($\gamma$) cannot occur for a non-singular pair.
\smallskip

The case ($\delta$): From $a=-a^{\sigma}(\zeta\eta/\lambda)$ and its $\sigma$-image $a^{\sigma}=-a(\zeta\eta/ \lambda)^{\sigma} = a(\zeta\eta/\lambda)^{-1}$ we have $a=-a$, that is, $a=0$. Similarly, from $d=d^{\sigma}\lambda/(\zeta\eta)$ one concludes that $d=0$, contrary to non-singularity.
\smallskip

We can now complete the enumeration of non-singular quadruples $(\lambda,a,d,u)$ in the case $B=\off(\lambda,1)$ for $\lambda=\xi^{(q+1)/2}$. As we saw above, the case $(\delta)$ is excluded. The case $(\gamma)$ does not yield  non-singular pairs but we still need more details about the corresponding quadruples to be excluded. We saw that $(\gamma)$ occurs if and only if $d=a\lambda/(\zeta\eta)$, and this gives $u=ad+\lambda^{\sigma}=a^2\xi^{(q-1)/2}\xi^{i(q-1)} + \lambda^{\sigma}$, $i\in \{0,1,\ldots,q\}$. Since now $q\equiv 1$ mod $4$, it follows that $u-\lambda^{\sigma}\in S(F)$, that is, $u\in N(\lambda^{\sigma})$. Conversely, if $u\in N(\lambda^{\sigma})$, then for each of the $q+1$ values of $i\in \{0,1,\ldots,q\}$ we obtain two solutions for $a$ as above (and for each such $a$ a unique $d$), not giving a non-singular pair. It follows that for each of the $n_F$ elements $u\in N(\lambda^{\sigma})$ we can choose $(q^2-1)-2(q+1)=(q+1)(q-3)$ values $a\in F^*$ (and for each such $a$ a unique $d$) such that $ad+\lambda^{\sigma}=u$, so that the order of $[A,1]$ is not equal to $4$. This gives a total of $n_F(q+1)(q-3)$ quadruples $(\xi^{(q+1)/2},a,d,u)$ collected so far. Any of the remaining $(q^2-3)/2 - n_F$ non-square values of $u$ distinct from $\lambda^{\sigma}$ then give $((q^2-3)/2 - n_F)(q^2-1)$ quadruples $(\xi^{(q+1)/2},a,d,u)$ that we have to consider for non-singular pairs. Finally, if $u=\lambda^{\sigma}$, then either $a=0$ with $q^2-1$ choices for $d$, or $d=0$ with $q^2-1$ choices for $a$, giving further $2(q^2-1)$ quadruples. Since we have shown earlier that the induced action of the group $H^*$ on these quadruples is semi-regular with orbits of length $4(q+1)$, the number $n_4$ of $\ovl{G}$-equivalence classes of non-singular quadruples $(\lambda,a,d,u)$ for off-diagonal $B$ and $\lambda=\xi^{(q+1)/2}$, $q\equiv 1$ mod $4$, is
\[ n_4=  \frac{n_F(q+1)(q-3) + ((q^2-3)/2 - n_F)(q^2-1)+2(q^2-1)}{4(q+1)} \ ,\]
which simplifies to
\begin{equation}\label{ad2} n_4=  \frac{1}{8}\left((q-1)(q^2+1) - 4n_F\right)\ .\end{equation}

%------------------------------------------------
\section{Enumeration of orientably-regular maps on {\bf $M(q^2)$}}\label{sec:enum}
\smallskip

We have seen in section \ref{sec:twpairs} that a pair $([A,1],[B,1])$ of elements of $G$, with diagonal $B$ and $A$ given by (\ref{A1}) or with off-diagonal $B$ and $A$ given by (\ref{A2}), and with product of order two,  generates a twisted subgroup of $G$ if and only if the pair is non-singular. In the previous two sections we have counted orbits of non-singular pairs in $G$ under conjugation in $\ovl{G}$, with no regard to subgroups the pairs generate. The number of these orbits turns out to be $n_1+n_3+n_4$ if $q\equiv 1$ mod $4$ and $n_1+n_2+n_3$ if $q\equiv -1$ mod $4$; in both cases the sum is equal to $(q^2-1)(q^2-2)/8$. We state this as a separate result.
\smallskip

\bp\label{count}
The number of $\ovl{G}$-orbits of non-singular pairs in $G=M(q^2)$ is equal to $(q^2-1)(q^2-2)/8$. \hfill $\Box$
\ep

We will now refine our considerations and take into account subgroups generated by non-singular pairs. For our group $G=G_{2f}=M(p^{2f})$ and for any positive divisor $e$ of $f$ such that $f/e$ is odd we let $G_{2e}$ denote the canonical copy of $M(p^{2e})$ in $G$. In Lemma \ref{subgps} we saw that the automorphism $\sigma$ of $F=F_{2f}=\GF(p^{2f})$ of order two restricts to an automorphism $\sigma_{2e}$ of order two of the subfield $F_{2e}=\GF(p^{2e})$ of $F$. We recall that $\ovl{G}=\ovl{G}_{2f}= G_{2f}\langle\sigma\rangle$ and we similarly introduce $\ovl{G}_{2e}$ for every $e$ as above by letting  $\ovl{G}_{2e}=G_{2e}\langle\sigma_{2e}\rangle$.
\smallskip

Let $\orb_f(e)$ denote the number of $\ovl{G}_{2f}$-orbits of non-singular pairs $([A,1],[B,1])$ of $G$ that generate a subgroup of $G$ isomorphic to $M(p^{2e})$. At the same time, let $\orb(e)$ be the number of orbits of non-singular pairs of $G_{2e}$ which generate $G_{2e}$. The two quantities, are, in fact, equal, which is fundamental for our enumeration.
\smallskip

\bp\label{equal}
For each positive divisor $e$ of $f$ with $f/e$ odd, we have $\orb_f(e)=\orb(e)$.
\ep

\pr
It is clear that every $\ovl{G}_{2e}$-orbit of a non-singular pair in the canonical copy $G_{2e}\cong M(p^{2e})$ in $G_{2f}$ is contained in a $\ovl{G}_{2f}$-orbit of the same pair. In the reverse direction, let a non-singular pair in $G$ generate a subgroup isomorphic to $M(p^{2e})$. Since, by the important Proposition \ref{e/f-conjug-c}, all such subgroups are $\ovl{G}_{2f}$-conjugate in $G_{2f}$, we may assume that the non-singular pair is contained in $G_{2e}$. But then the $\ovl{G}_{2f}$-orbit of this pair obviously contains a $\ovl{G}_{2e}$-orbit of the same pair. The result will now be a consequence of the following claim:
\smallskip

{\sl
Let $([A,1],[B,1])$ and $([A',1],B',1])$ be two non-singular pairs of $G_{2e}$ both generating $G_{2e}$ and lying in the same $\ovl{G}_{2f}$-orbit of $G_{2f}$. Then the two pairs are contained in the same $\ovl{G}_{2e}$-orbit of $G_{2e}$. }
\smallskip

We have $B=\dia(\lambda,1)$ and $B'=\dia(\lambda',1)$, or $B=\off(\lambda,1)$ and $B'=\off(\lambda',1)$, in both cases for some $\lambda,\lambda'\in F_{2e}$. The assumption that the two pairs are contained in the same $\ovl{G}_{2f}$-orbit means, by Proposition \ref{conj}, that there exists an element $P\in \GL(2,p^{2f})$ such that $AP^{\sigma}=\veps PA'^{(\sigma)}$ and $BP^{\sigma}=\delta PB'^{(\sigma)}$ for some $\veps,\delta\in F^*$. Propositions \ref{p-dia} and \ref{p-off} applied to the equation $BP^{\sigma}=\delta PB'^{(\sigma)}$ tell us that $P=\dia(\omega,1)$ or $P=\off(\omega,1)$ for some $\omega\in F^*$. We will show that for all choices of $A$, $A'$ and $P$ the equation $AP^{\sigma}=\veps PA'^{(\sigma)}$ implies that $\omega\in F_{2e}$, which will prove that our two non-singular pairs are in one $G_{2e}$-orbit.
\smallskip

There are two possibilities for $A$ and $A'$, given in (\ref{A1}) and (\ref{A2}), and two possibilities for $P$, the diagonal and the off-diagonal case, to substitute in the equation $AP^{\sigma}=\veps PA'^{(\sigma)}$; in all the four cases we need to conclude that $\omega\in F_{2e}$. We demonstrate it on the case when $P=\off(\omega,1)$ and when $A, A'$ are given by (\ref{A2}), that is, when $A$ has rows $(a,\lambda^{\sigma})$, $(-1,d)$ and $A'$ has rows $(a',\lambda'^{\sigma})$, $(-1,d')$ for some $a,d,\lambda,a',d',\lambda'\in F_{2e}$. The above equation then reads
$$
\begin{pmatrix} a & \lambda^{\sigma} \\ -1 & d\end{pmatrix}
\begin{pmatrix} 0 & \omega^{\sigma} \\ 1 & 0 \end{pmatrix} = \veps \begin{pmatrix} 0 & \omega \\ 1 & 0 \end{pmatrix}
{\begin{pmatrix} a' & \lambda'^{\sigma} \\ -1 & d'\end{pmatrix}}^{(\sigma)} \ , $$
giving the equations
$$ \lambda^{\sigma}=-\veps \omega\ , \ \ \ \omega^{\sigma}a=\veps\omega d'^{(\sigma)}\ , \ \ \
d = \veps a'^{(\sigma)}\ , \ \ \ -\omega = \veps^{\sigma}\lambda'^{(\sigma)}\ . $$
If $a', d\ne 0$, then the third equation in combination with the first gives $\omega = -\lambda^{\sigma}/\veps = -\lambda^{\sigma}a'^{(\sigma)}/d$, showing that $\omega\in F_{2e}$. If $a'=d=0$, we must have $a,d'\ne 0$ and then a combination of the first two equations gives $\omega^{\sigma}= \veps\omega d'^{(\sigma)}/a = -\lambda^{\sigma}d'^{(\sigma)}/a$, showing again that $\omega\in F_{2e}$. The remaining three cases are similar (and easier) to check. This completes the proof.
\ebox

For positive integers $x$ let us define a function $h$ by $h(x)= (p^{2x}-1)(p^{2x}-2)/8$. In terms of $h$ and the numbers $\orb_f(e)$, Proposition \ref{count} simply says that $\sum_{e}\orb_f(e)=h(f)$, where summation is taken over all positive divisors $e$ of $f$ such that $f/e$ is odd. By Proposition \ref{equal} we may replace $\orb_f(e)$ with $\orb(e)$ and obtain $\sum_{e}\orb(e)=h(f)$, with the same summation convention. This miniature but important detail enables us to make a substantial advance in the enumeration.
\smallskip

Let $f=2^{\alpha}o$ where $o$ is an odd integer and let $e$ be a divisor of $f$ such that $f/e$ is odd; equivalently, $e=2^{\alpha}d$ where $d$ is a positive (and necessarily odd) divisor of $o$. Taking the above notes into account, Proposition \ref{count} may then be restated as follows:
\be\label{eq_m} \sum_{d\mid o} \orb (2^{\alpha}d) = h(2^{\alpha}o)\ .\ee
Using the M\"obius inversion we obtain $\orb(f)= \orb(2^{\alpha}o) =\sum_{d\mid o} \mu(o{/}d)h(2^{\alpha}{}d)$, where $\mu$ is the classical number-theoretic M\"obius function $\mu$ on positive integers. We thus arrive at our first main result.
\smallskip

\bt\label{main1} Let $q=p^{f}$ for an odd prime $p$, let $G=M(q^2)$, and let $f=2^{\alpha}{}o$ with $o$ odd. The number of $\ovl{G}$-orbits of non-singular generating pairs of $G$ is equal to
\[\sum_{d\mid o} \mu(o{/}d)h(2^{\alpha}{}d)\ , \ \  {\rm where} \ \  h(x)= (p^{2x}-1)(p^{2x}-2)/8\ .\]
\et

The last step is to study conjugacy of non-singular pairs of $M(q^2)$ under the action of the group ${\rm Aut}(M(q^2))$ which, as we know by Proposition \ref{cor-auto}, is isomorphic to ${\rm P\Gamma L}(2,q^2)$. Since for $q=p^{f}$ we have ${\rm P\Gamma L}(2,q^2)\cong {\rm PGL}(2,q^2)\rtimes Z_{2f}\cong \ovl{G}\rtimes Z_f$, it is sufficient to investigate the induced action of the Galois automorphisms $\sigma_j:\ z\mapsto z^{p^j}$ for $z\in F=\GF(p^{2f})$ and $1\le j\le f-1$ on the $\ovl{G}$-orbits of our non-singular pairs $([A,1],[B,1])$. We will use the natural notation $O^{\sigma_j}$ for the $\sigma_j$-image of a $\ovl{G}$-orbit $O$ of a pair $([A,1],[B,1])$ of elements of $G$. Note that $\sigma_{f}=\sigma$, and we also have $O^{\sigma_f}=O$, by the remark made at the beginning of section \ref{sec:twpairs}. Clearly, if $O^{\sigma_j}\cap O\ne \emptyset$, then $O^{\sigma_{j}}=O$.

\bp\label{Galois}
Let $O$ be the orbit of a non-singular pair $([A,1],[B,1])$ of elements of $G$ under conjugation in $\ovl{G}$ and let $j$ be the smallest positive integer for which $O^{\sigma_j}=O$. If $[A,1]$ and $[B,1]$ generate $G$, then $j=f$.
\ep

\pr
We may assume that $f\ge 2$, otherwise the result is trivial. Suppose that the pair $([A,1],[B,1])^{\sigma_j}= ([A,1]^{\sigma_j}, [B,1]^{\sigma_j})$ is $\ovl{G}$-conjugate to the pair $([A,1],[B,1])$, that is, there exists some $C\in {\rm GL}(2,q^2)$ and $i\in Z_2$ such that  $[A,1]^{\sigma_j}=[C,i]^{-1}[A,1][C,i]$ and $[B,1]^{\sigma_j}= [C,i]^{-1}[B,1][C,i]$. It follows that for every $[X,1]\in \langle [A,1],[B,1] \rangle$ we have $[X,1]^{\sigma_j}=[C,i]^{-1}[X,1][C,i]$. Using our assumption that $\langle [A,1],[B,1] \rangle = G$, we conclude that the above is valid also for $X=\dia(\xi,1)$, where $\xi$ is a primitive element of $F=\GF(p^{2f})$. Letting $C$ have elements $\alpha,\beta,\gamma,\delta$ in the usual order, the equivalent form $[C,i][\dia(\xi^{p^j},1),1]=[\dia(\xi,1),1][C,i]$ of the above equation yields
$$ \begin{pmatrix} \alpha & \beta \\ \gamma & \delta\end{pmatrix}
\begin{pmatrix}  (\xi^{p^j})^{(\sigma)}& 0 \\ 0 & 1 \end{pmatrix} = \veps \begin{pmatrix} \xi & 0 \\ 0 & 1 \end{pmatrix}
\begin{pmatrix} \alpha^{\sigma} & \beta^{\sigma} \\ \gamma^{\sigma} & \delta^{\sigma}\end{pmatrix}  $$
for some $\veps \in F^*$; here we used the $(\sigma)$-convention introduced at the end of section \ref{sec:diag}.    This gives the system of equations
$$ \alpha (\xi^{p^j})^{(\sigma)}=\veps\xi\alpha^{\sigma} \ , \ \   \beta=\veps\xi\beta^{\sigma} \ , \ \
\gamma (\xi^{p^j})^{(\sigma)}=\veps\gamma^{\sigma} \ , \delta=\veps\delta^{\sigma} \ .$$
Consider first the case $\delta\ne 0$; without loss of generality we then may assume $\delta=1$. Then $\veps=1$, and the equation for $\alpha$ gives $(\xi^{p^j})^{(\sigma)}\xi^{-1} = \alpha^{\sigma}\alpha^{-1}$, or, equivalently, $\xi^{p^{j+if}-1}=\alpha^{p^f-1}$. It follows that $p^f-1$ is a divisor of $p^{j+if}-1$, which implies that $f$ divides $j+if$ and hence $f$ divides $j$, which, since $j\le f$, shows that $j=f$. If $\delta=0$ then, without loss of generality, $\beta=1$ and so $\veps=\xi^{-1}$. The equation for $\gamma$ now implies $\xi^{p^{j+if}+1}=\gamma^{p^f-1}$. It follows that $p^f-1$ divides $p^{j+if}+1$ and hence also $p^{2(j+if)}-1$ and so $f$ must divide $2(j+if)$. Thus, $f$ is a divisor of $2j$ and as $j\le f$, we have either $j=f$ or $j=f/2$ (assuming $f$ is even). But the last case is easily seen to be impossible since $p^f-1$ is not a divisor of $p^{f/2}+1$ or $p^{3f/2}+1$. This completes the proof.
\ebox

Proposition \ref{Galois} tells us that if a non-singular pair $([A,1],[B,1])$ of elements of $G$ actually generates $G$ and gives rise to an orbit $O$ under conjugation in $\ovl{G}$, then the action of the group ${\rm Aut}(M(q^2))$ fuses the $f$ orbits $O^{\sigma_j}$ for $j\in \{0,1,\ldots,f-1\}$ into a single orbit. Recalling the one-to-one correspondence between isomorphism classes of orientably-regular maps supported by the group $G=M(q^2)$ and orbits of (necessarily non-singular) generating pairs of $G$ under conjugation by ${\rm Aut}(G)$, Theorem \ref{main1} then immediately implies our second main result.

\bt\label{main2} Let $q=p^{f}$ for an odd prime $p$ and let $f=2^{\alpha}{}o$ with $o$ odd. The number of orbits of non-singular generating pairs of $M(q^2)$ under the action of the group ${\rm Aut}(M(q^2))$, and hence the number of isomorphism classes of orientably-regular maps ${\cal M}$ with ${\rm Aut}^+({\cal M})\cong M(q^2)$, is equal to
\[ \frac{1}{f}\sum_{d\mid o} \mu(o{/}d)h(2^{\alpha}{}d)\ , \ \  {\rm where} \ \  h(x)= (p^{2x}-1)(p^{2x}-2)/8\ .\]
\et

%------------------------------------------------
\section{Enumeration of reflexible maps}\label{sec:ref}

A map is called \emph{reflexible} if it admits an automorphism reversing the orientation of the surface. For orientably-regular maps represented by triples $(G,x,y)$ as indicated in the Introduction, reflexibility is equivalent with the existence of an automorphism $\theta$ of the group $G$ such that $\theta(x)=x^{-1}$ and $\theta(y)=y^{-1}$. Note that if such a $\theta$ exists, then $\theta^2=id$.
\smallskip

In the specific situation considered in this paper, namely, when $G=M(q^2)$ for $q=p^f$, we established in Proposition \ref{cor-auto} that ${\rm Aut}(G)\cong {\rm P\Gamma L}(2,q^2)$. Moreover, it is well known that every automorphism in ${\rm P\Gamma L}(2,q^2) \cong {\rm PGL}(2,q^2)\rtimes Z_{2f}$ is a composition of a conjugation by some element of ${\rm PGL}(2,q^2)$ and a power of the Frobenius automorphism $z\mapsto z^p$ of the Galois field $F=GF(p^{2f})$. It follows that an {\sl involutory} automorphism $\theta$ of $G=M(q^2)$ is a composition of a conjugation as above with $\sigma^i$ for $i\in \{0,1\}$, where $\sigma$ is the automorphism of $F$ sending $z$ to $z^q$. By the remark at the beginning of Section \ref{sec:twpairs}, however, the action of $\sigma$ is equivalent to conjugation in $\ovl{G}=G\langle\sigma\rangle$ by the element $[I,1]$. Consequently, an orientably-regular map on the group $G=M(q^2)$ generated by a pair of elements $x=[A,1]$ and $y=[B,1]$ is reflexible if and only if the ordered pairs $(x,y)$ and $(x^{-1},y^{-1})$ are conjugate by an {\sl involutory} element of $\ovl{G}$.
\smallskip

In this section we will count the number of reflexible orientably-regular maps on $M(q^2)$. In particular, we will see that not all orientably-regular maps with automorphism group $M(q^2)$ are reflexible, in contrast with the position for $PGL(2,q^2)$, see e.g. \cite{CPS2}. We will use techniques similar to the main enumeration in previous sections and structure our explanations accordingly.

\subsection{Conjugating involutions}

In view of the above findings, an orientably-regular map on $M(q^2)$ given by the generating pair $([A,1],[B,1])$ is reflexible if and only if there is an involution $[C,i]\in\overline{G}$ for some $i\in\{0,1\}$ such that $[C,i][A,1][C,i]=[A,1]^{-1}$ and $[C,i][B,1][C,i]=[B,1]^{-1}$. As before, we will be dealing separately with the diagonal and off-diagonal cases for the matrix $B$.

\subsubsection{The diagonal case $B=\dia(\lambda,1)$}

We can write $\left(B^{\sigma}\right)^{-1}=\begin{pmatrix}\veps \lambda^{-q}&0\\0&\veps\end{pmatrix}$ for some $\veps\neq 0$. We seek a conjugating involution $[C,i]$ as above and begin with $i=0$. For $[C,0]$ to be an involution we must have $\mathrm{tr}(C)=0$ and we can without loss of generality write
$C=\begin{pmatrix}1&\beta\\ \gamma &-1\end{pmatrix}$ or $C=\begin{pmatrix}0&\beta\\1&0\end{pmatrix}$.
\smallskip

For the first form of $C$ the conjugation equation becomes:
\[
\begin{pmatrix}1&\beta\\ \gamma &-1\end{pmatrix}
\begin{pmatrix}\lambda&0\\0&1\end{pmatrix}
\begin{pmatrix}1&\beta^q\\ \gamma^q&-1\end{pmatrix}
=\begin{pmatrix}\veps \lambda^{-q}&0\\0&\veps\end{pmatrix}
\]
This leads to the system of equations $\lambda+\beta\gamma^q=\veps\lambda^{-q}$, $\lambda \beta^q-\beta=0$, $\gamma\lambda-\gamma^q=0$, and $\beta^q \gamma\lambda+1=\veps$. If $\beta=0$ then $\veps=1$ and so $\lambda=\lambda^{-q}$ or $\lambda\in\sqrt[q+1]{1}$. If $\beta\neq 0$ then $\lambda=\beta^{1-q}$. In either case, $\lambda\in S(F)$ contrary to the construction of $B$. Hence there are no conjugating matrices $C$ of this form.
\smallskip

For the second possibility we have:
\[
\begin{pmatrix}0&\beta\\1&0\end{pmatrix}
\begin{pmatrix}\lambda&0\\0&1\end{pmatrix}
\begin{pmatrix}0&\beta^q\\1&0\end{pmatrix}
=\begin{pmatrix}\veps \lambda^{-q}&0\\0&\veps\end{pmatrix}
\]
This leads to $\beta=\veps\lambda^{-q}$ and $\lambda \beta^q=\veps$, giving potential solutions of this form with $\beta$ satisfying $\beta^{q-1}=\lambda^{q-1}$.
\smallskip

Now consider the case $i=1$. Since $[C,1]$ is an involution, it cannot be in $G$; it follows that  $CC^{\sigma}=\veps I$ for some $\veps\neq 0$ and $\mathrm{det}(C)\in S(F)$. Letting $C=\begin{pmatrix}\alpha &\beta \\ \gamma & \delta\end{pmatrix}$ we may either have $\gamma=0$ in which case we can take $\alpha=1$, or $\gamma=1$, furnishing again two possible forms for $C$.
\smallskip

The first form of $C$ for $i=1$ we will consider is $C=\begin{pmatrix}1&\beta\\0&\delta\end{pmatrix}$ for some $\delta\in S(F)$. Since $CC^{\sigma}=\veps I$ we must have $\veps=1$ and $\delta\in\sqrt[q+1]{1}$. In this case the matrix equations are:
\[
\begin{pmatrix}1&\beta\\0&\delta\end{pmatrix}
\begin{pmatrix}\lambda^q&0\\0&1\end{pmatrix}
\begin{pmatrix}1&\beta\\0&\delta\end{pmatrix}
=\begin{pmatrix}\veps \lambda^{-q}&0\\0&\veps\end{pmatrix}
\]
It follows that $\lambda^q=\veps\lambda^{-q}$ and $\delta^2=\veps$, implying $\delta^2=\lambda^{2q}$, but as $\delta\in S(F)$ and $\lambda\in N(F)$ this is impossible. We conclude that the first form cannot occur.

The second form of $C$ for $i=1$ to consider is $C=\begin{pmatrix}\alpha&\beta\\1&\delta\end{pmatrix}$ with $\alpha\delta-\beta\in S(F)$. The requirement for $[C,1]$ to have order two gives:
\[
\begin{pmatrix}\alpha&\beta\\1&\delta\end{pmatrix}
\begin{pmatrix}\alpha^q&\beta^q\\1&\delta^q\end{pmatrix}
=\begin{pmatrix}\veps&0\\0&\veps\end{pmatrix}
\]
Thus, $\delta=-\alpha^q$ and $C$ has the form $\begin{pmatrix}\alpha&\beta\\1&-\alpha^q\end{pmatrix}$. The conjugation equation then becomes:
\[
\begin{pmatrix}\alpha&\beta\\1&-\alpha^q\end{pmatrix}
\begin{pmatrix}\lambda^q&0\\0&1\end{pmatrix}
\begin{pmatrix}\alpha&\beta\\1&-\alpha^q\end{pmatrix}
=\begin{pmatrix}\veps \lambda^{-q}&0\\0&\veps\end{pmatrix}
\]
Now,$\alpha\lambda^q-\alpha^q=0$. If $\alpha\neq 0$, then $\lambda^q=\alpha^{q-1}\in S(F)$, which is impossible. So $\alpha=0$ and $C$ must have the form $\begin{pmatrix}0&\beta\\1&0\end{pmatrix}$ for some $\beta$. Since $CC^{\sigma}=\veps I$ we must have $\beta=\beta^q$ and hence $\beta\in F'^*$.

\subsubsection{The off-diagonal case $B=\mathrm{off}(\lambda,1)$}

We can write $\left(B^{\sigma}\right)^{-1}=\begin{pmatrix}0&\veps\\ \veps\lambda^{-q}&0\end{pmatrix}$ for some $\veps\neq 0$. We begin with the case $i=0$, and as before we can without loss of generality set
$C=\begin{pmatrix}1&\beta\\ \gamma&-1\end{pmatrix}$ or $C=\begin{pmatrix}0&\beta\\1&0\end{pmatrix}$.
For the first form of $C$ we have:
\[
\begin{pmatrix}1&\beta\\ \gamma&-1\end{pmatrix}
\begin{pmatrix}0&\lambda\\1&0\end{pmatrix}
\begin{pmatrix}1&\beta^q\\ \gamma^q&-1\end{pmatrix}
=\begin{pmatrix}0&\veps\\ \veps\lambda^{-q}&0\end{pmatrix}
\]
It follows that  $\beta+\gamma^q \lambda=0$, $\beta^{q+1}-\lambda=\veps$, $\gamma^{q+1}\lambda-1=\veps\lambda^{-q}$ and $-\beta^q-\gamma\lambda=0$. If $\beta=0$ then $\gamma=0$ and so $\veps=-\lambda$ and $\lambda^{1-q}=1$. If $\beta,\gamma\neq 0$, then $\beta+\gamma^q\lambda^q=0$, giving $\lambda=\lambda^q$. In either case we conclude that $\lambda\in S(F)$, contrary to the construction of $B$. Consequently, there are no matrices $C$ of this form.
\smallskip

For the possibility $C=\begin{pmatrix}0&\beta\\1&0\end{pmatrix}$ we have:
\[
\begin{pmatrix}0&\beta\\1&0\end{pmatrix}
\begin{pmatrix}0&\lambda\\1&0\end{pmatrix}
\begin{pmatrix}0&\beta^q\\1&0\end{pmatrix}
=\begin{pmatrix}0&\veps\\ \veps\lambda^{-q}&0\end{pmatrix}
\]
This leads to $\beta^{q+1}=\veps$, $\lambda=\veps\lambda^{-q}$ and we have potential solutions of this form with $\beta$ satisfying $\beta^{q+1}=\lambda^{q+1}$.
\smallskip

Now consider the second case $i=1$. As before the first possible form for $C$ is $\begin{pmatrix}1&\beta\\0&\delta\end{pmatrix}$ for some $\delta\in \sqrt[q+1]{1}$.
In this case the matrix equations are:
\[
\begin{pmatrix}1&\beta\\0&\delta\end{pmatrix}
\begin{pmatrix}0&\lambda^q\\1&0\end{pmatrix}
\begin{pmatrix}1&\beta\\0&\delta\end{pmatrix}
=\begin{pmatrix}0&\veps\\ \veps\lambda^{-q}&0\end{pmatrix}
\]
We must have $\beta=0$ and $\delta=\veps\lambda^{-q}$, giving potential solutions of this form for any $\delta\in \sqrt[q+1]{1}$.
\smallskip

The second possible form of $C$ for $i=1$ is $C=\begin{pmatrix}\alpha&\beta\\1&-\alpha^q\end{pmatrix}$ for which the equations are:
\[
\begin{pmatrix}\alpha&\beta\\1&-\alpha^q\end{pmatrix}
\begin{pmatrix}0&\lambda^q\\1&0\end{pmatrix}
\begin{pmatrix}\alpha&\beta\\1&-\alpha^q\end{pmatrix}
=\begin{pmatrix}0&\veps\\ \veps\lambda^{-q}&0\end{pmatrix}
\]
Here we have $\alpha\beta+\alpha\lambda^q=0$, $\beta^2-\alpha^{q+1}\lambda^q=\veps$, $-\alpha^{q+1}+\lambda^q= \veps\lambda^{-q}$, and $-\alpha^q\beta-\alpha^q\lambda^q=0$. If $\alpha\neq 0$, then $\lambda^q=-\beta$ but $\beta\in F'^*$, so $-\beta\in S(F)$ and this is impossible. It follows that $\alpha=0$ but then $\veps=\beta^2 =\lambda^{2q}$, while at the same time $\beta\in S(F)$ and $\lambda^q\in N(F)$, a contradiction. Hence there are no matrices $C$ of this form.

\subsubsection{Summary of possible conjugating elements}

In the case $B=\dia(\lambda,1)$ the possible conjugating elements can have the following forms.
\begin{equation}[C,0];\quad C=\begin{pmatrix}0&\beta\\1&0\end{pmatrix};\quad \beta^{q-1}=\lambda^{q-1}\end{equation}
\begin{equation}[C,1];\quad C=\begin{pmatrix}0&\beta\\1&0\end{pmatrix};\quad \beta\in F'^*\end{equation}
In the case $B=\off(\lambda,1)$ the possible conjugating elements can have the following forms.
\begin{equation}[C,0];\quad C=\begin{pmatrix}0&\beta\\1&0\end{pmatrix};\quad \beta^{q+1}=\lambda^{q+1}\end{equation}
\begin{equation}[C,1];\quad C=\begin{pmatrix}1&0\\0&\delta\end{pmatrix};\quad \delta\in \sqrt[q+1]{1}\end{equation}

\subsection{Enumeration}

We now proceed to the actual enumeration and follow the same strategy as we used for general maps, namely, counting orbits in the diagonal and off-diagonal cases for $B$ as in Sections \ref{sec:rep-dia} and \ref{sec:rep-off} and then deriving the final enumeration result using the M\"obius inversion formula as in Section \ref{sec:enum}.

\subsubsection{Counting orbits: The case $B=\dia(\lambda,1)$}

In this case we have $A=\begin{pmatrix}-1&b\\c&\lambda^q\end{pmatrix}$. We can write $[A,1]^{-1}$ as $[(A^{\sigma})^{-1},1]$ where $(A^{\sigma})^{-1}$ has the form $\veps\begin{pmatrix}-\lambda&b^q\\ c^q&1\end{pmatrix}$ for some $\veps\neq 0$.
\smallskip

The first possibility is that the conjugating element has the form $[C,0]$ with $C=\begin{pmatrix}0&\beta\\ 1&0\end{pmatrix}$ for some $\beta$ such that $\beta^{q-1}=\lambda^{q-1}$. In this case the equations are:
\[
\begin{pmatrix}0&\beta\\1&0\end{pmatrix}
\begin{pmatrix}-1&b\\c&\lambda^q\end{pmatrix}
\begin{pmatrix}0&\beta^q\\1&0\end{pmatrix}
=\veps\begin{pmatrix}-\lambda&b^q\\c^q&1\end{pmatrix}
\]
After some manipulation the resulting equations lead to $\beta=-b^q/c;\,\beta^{q-1}=\lambda^{q-1}$.
Another way to express the second equation is $(bc)^{1-q}=\lambda^{q-1}$ or $(bc\lambda)^{q-1}=1$ so that $bc\lambda\in F'^*$.
\smallskip

Since $u=bc+\lambda^q\in N(F)\setminus\{\lambda^q\}$ we can count the possible values of $u$ by noting that $u\lambda=bc\lambda+\lambda^{q+1}$. So $u\lambda\in F'^*$ but we must have $u\neq\lambda^q$ and therefore there are exactly $q-2$ choices for $u$. For each choice of $u$ we then set $bc=u-\lambda^q$ and then there are $q^2-1$ choices for $b$ and $c$ is determined.
\smallskip

For given $\lambda=\xi^i$, the number of matrices $A$ satisfying the equations is then $(q-2)(q^2-1)$.
The $\ovl{G}$-orbits have length $2(q-1)$ if $i<(q-1)/2$ and $4(q-1)$ if $i=(q-1)/2$.
So the number of orbits of generating pairs of this form for each $i$ is:
\begin{equation}
r_{1,i}=
\begin{cases}
\displaystyle (q-2)\frac{q+1}{2}&\quad\text{ for }\lambda=\xi^i,i<(q-1)/2\\[1em]
\displaystyle (q-2)\frac{q+1}{4}&\quad\text{ for }\lambda=\xi^i,i=(q-1)/2
\end{cases}
\label{eq:ref1a}
\end{equation}
The total number of orbits is therefore
\begin{equation}
r_1=\sum_{i\leq(q-1)/2}r_{1,i}=\frac{(q^2-1)(q-2)}{8}
\label{eq:ref1}
\end{equation}

The second possibility is that the conjugating element has the form $[C,1]$ with $C=\begin{pmatrix}0&\beta\\ 1&0\end{pmatrix}$ for some $\beta\in F'^*$. In this case the equations are:
\[
\begin{pmatrix}0&\beta\\1&0\end{pmatrix}
\begin{pmatrix}-1&b^q\\c^q&\lambda\end{pmatrix}
\begin{pmatrix}0&\beta\\1&0\end{pmatrix}
=\veps\begin{pmatrix}-\lambda&b^q\\c^q&1\end{pmatrix}
\]
The resulting equations become $\beta=-b^q/c^q;\,\beta\in F'^*$; in particular, $b/c\in F'^*$. Now note that $bc=c^2b/c\in S(F)$ and so the number of possible values for $u=bc+\lambda^q$ is just $n_F=(q^2-1)/4$ from before.
It follows that there are $(q^2-1)/4$ possible values for the product $bc$ and we can count the number of ways we can choose $b$.
\smallskip

Since $b/c\in F'^*$ we can write $x=b/c,\ y=bc$ and there are $q-1$ choices for $x$. For each choice we need $b$ to satisfy the equations $x=b/c,\ y=bc$ so $xy=b^2$. Knowing that $xy\in S(F)$, there are exactly two choices for $b$ and so the total number of matrices $A$ satisfying the equations for a given $\lambda$ is $(q-1)(q^2-1)/2$.
\smallskip

Then, given the orbit lengths as before, the number of orbits of pairs of this form for each $i$ is:
\begin{equation}
r_{2,i}=
\begin{cases}
\displaystyle \frac{q^2-1}{4}&\quad\text{ for }\lambda=\xi^i,i<(q-1)/2\\[1em]
\displaystyle \frac{q^2-1}{8}&\quad\text{ for }\lambda=\xi^i,i=(q-1)/2
\end{cases}
\label{eq:ref2a}
\end{equation}
\begin{equation}
r_2=\sum_{i\leq(q-1)/2}r_{2,i}=\frac{(q^2-1)(q-1)}{16}
\label{eq:ref2}
\end{equation}
Note also that in the set of matrices counted by $r_1$ above, we have $bc\in N(F)$ and for those counted by $r_2$ we have $bc\in S(F)$, so that the sets are disjoint.

\subsubsection{Counting orbits: The case $B=\mathrm{off}(\lambda,1)$}

In this case we have $A=\begin{pmatrix}a&\lambda^q\\-1&d\end{pmatrix}$. We can write $[A,1]^{-1}$ as $[(A^{\sigma})^{-1},1]$ where $(A^{\sigma})^{-1}$ has the form $\veps\begin{pmatrix}d^q&-\lambda\\ 1&a^q\end{pmatrix}$ for some $\veps\neq 0$.
\smallskip

The first possibility is that the conjugating element has the form $[C,0]$ with $C=\begin{pmatrix}0&\beta\\ 1&0\end{pmatrix}$ for some $\beta$ such that $\beta^{q+1}=\lambda^{q+1}$. In this case the equations are:
\[
\begin{pmatrix}0&\beta\\1&0\end{pmatrix}
\begin{pmatrix}a&\lambda^q\\-1&d\end{pmatrix}
\begin{pmatrix}0&\beta^q\\1&0\end{pmatrix}
=\veps\begin{pmatrix}d^q&-\lambda\\1&a^q\end{pmatrix}
\]
It follows that $\beta d=\veps d^q$, $-\beta^{q+1}=-\veps\lambda$, $\lambda^q=\veps$, and $a\beta^q=\veps a^q$.
So either $ad=0$ or $(\lambda ad)^{q-1}=1$. If $ad\neq 0$ then as before there are exactly $q-2$ choices for $u=ad+\lambda^q$, each with $q^2-1$ choices for $a$ and then $d$ is determined. If $ad=0$ there are $q^2-1$ choices for $a$ if $d=0$ and $q^2-1$ choices for $d$ if $a=0$.
\smallskip

Thus for a given $\lambda=\xi^i$ there are $q(q^2-1)$ matrices $A$ satisfying the equations. The orbit lengths are $2(q+1)$ if $i<(q+1)/2$ and $4(q+1)$ if $i=(q+1)/2$. So the number of orbits of this form for each $i$ is:
\begin{equation}
r_{3,i}=
\begin{cases}
\displaystyle \frac{q(q-1)}{2}&\quad\text{ for }\lambda=\xi^i,i<(q+1)/2\\[1em]
\displaystyle \frac{q(q-1)}{4}&\quad\text{ for }\lambda=\xi^i,i=(q+1)/2
\end{cases}
\label{eq:ref3a}
\end{equation}
\begin{equation}
r_3=\sum_{i\leq(q+1)/2}r_{3,i}=\frac{q(q^2-1)}{8}
\label{eq:ref3}
\end{equation}

The remaining possibility is that the conjugating element has the form $[C,1]$ with $C=\begin{pmatrix}1&0\\0&\delta\end{pmatrix}$ for
some $\delta$ such that $\delta^{q+1}=1$. In this case the equations are:
\[
\begin{pmatrix}1&0\\0&\delta\end{pmatrix}
\begin{pmatrix}a^q&\lambda\\-1&d^q\end{pmatrix}
\begin{pmatrix}1&0\\0&\delta\end{pmatrix}
=\veps\begin{pmatrix}d^q&-\lambda\\1&a^q\end{pmatrix}
\]
This leads to $a^q=\veps d^q$, $\lambda\delta=-\lambda\veps$, $-\delta=\veps$, and $\delta^2 d^q=\veps a^q$.
So $\delta=-a^q/d^q$ and $\delta\in\sqrt[q+1]{1}$. Now $ad\in S(F)$ since $a^q/d^q=-\delta\in S(F)$, which implies that the number of possible $u=ad+\lambda^q$ is, as before, $n_F=(q^2-1)/4$. To proceed, write $x=a/d,\ y=ad$ and then $x^{q(q+1)}=1$. This gives $x\in\sqrt[q+1]{1}$ and so there are $q+1$ choices for $x$ and for each $x$ there are $(q^2-1)/4$ choices for $y$, with still having to choose $d$ to satisfy $d^2=y/x$. Since we know $y/x\in S(F)$ there are 2 choices for $d$.
\smallskip

Thus, for a given $\lambda=\xi^i$, there are $(q+1)(q^2-1)/2$ matrices $A$ satisfying the above equations.
The orbit lengths are as before, so the number of orbits of this form for each $i$ is:
\begin{equation}
r_{4,i}=
\begin{cases}
\displaystyle \frac{q^2-1}{4}&\quad\text{ for }\lambda=\xi^i,i<(q+1)/2\\[1em]
\displaystyle \frac{q^2-1}{8}&\quad\text{ for }\lambda=\xi^i,i=(q+1)/2
\end{cases}
\label{eq:ref4a}
\end{equation}
\begin{equation}
r_4=\sum_{i\leq(q+1)/2}r_{4,i}=\frac{(q+1)(q^2-1)}{16}
\label{eq:ref4}
\end{equation}
Again, the sets of matrices counted by $r_3$ and $r_4$ are disjoint since in the first case $ad\in N(F)$ and in the second $ad\in S(F)$.

\subsubsection{Summary of counting orbits}\label{subsec:orb}

For $B=\mathrm{dia}(\lambda,1)$ the total number of $\ovl{G}$-orbits of generating pairs for reflexible maps is
\begin{equation}
R_1=r_1+r_2=\frac{(q^2-1)(3q-5)}{16}
\label{eq:totdia}
\end{equation}
For $B=\mathrm{off}(\lambda,1)$ the total number of such orbits is
\begin{equation}
R_2=r_3+r_4=\frac{(q^2-1)(3q+1)}{16}
\label{eq:totoff}
\end{equation}
The total number of orbits is therefore:
\begin{equation}
R=R_1+R_2=\frac{(q^2-1)(3q-2)}{8}
\label{eq:tot}
\end{equation}

\subsection{Counting reflexible maps}

We may enumerate the orientably-regular reflexible maps on $M(q^2)$ by using the above calculations in conjunction with the logic of Section~\ref{sec:enum}. Since details of this process are exactly as in Section~\ref{sec:enum} except for using the input on counting orbits from Subsection \ref{subsec:orb} we present just the final result.

\begin{thm}\label{main3}
Let $q=p^f$ be an odd prime power, with $f=2^{\alpha}o$ where $o$ is odd.
The number of orientably-regular reflexible maps $\cal M$ with $\Aut^+({\cal M})\cong M(q^2)$ is, up to isomorphism, equal to
\[ \frac{1}{f}\sum_{d\mid o} \mu(o{/}d)\tilde h(2^{\alpha}{}d)\ ,\]
where $\tilde h(x)= (p^{2x}-1)(3p^{x}-2)/8$ and $\mu$ is the M\"obius function.
\end{thm}

%------------------------------------------------
\section{Remarks}\label{sec:rem}

As stated in the Introduction, orientably-regular maps have been enumerated for very sparse classes of non-trivial groups, and in terms of {\sl all} such maps (not just maps of restricted types) only for almost-Sylow-cyclic groups \cite{CPS} and the linear fractional groups ${\rm PSL}(2,q)$ and ${\rm PGL}(2,q)$ \cite{CPS2,Sah}. It should be noted, however, that the available results for ${\rm PSL}(2,q)$ and ${\rm PGL}(2,q)$ are more detailed by giving `closed formulae' for the number of orientably-regular maps of every given type, whereas our main results in Theorems \ref{main2} and \ref{main3} contain formulae for the total number of such maps.
\smallskip

In order to obtain a refined version of our enumeration of orientably-regular maps with automorphism group isomorphic to a twisted linear fractional group $G=M(q^2)$ one could follow \cite{Jon}, which requires setting up both a character table for $G$ and the M\"obius function for the lattice of subgroups of $G$. The number of orientably-regular maps on the group $G$ is then obtained as a combination of a character-theoretic formula for counting solutions of the equation $xyz=1$ for $x,y,z$ in given conjugacy classes of $G$ (a special case of a general formula of Frobenius \cite{Fro}) combined with M\"obius inversion, which is a forthcoming project of the authors. Whether the project will return a `nice' formula, however, is not clear due to another significant difference between the family of orientably-regular maps on $M(q^2)$ compared to those on ${\rm PGL}(2,q)$. Namely, in the case of ${\rm PGL}(2,q)$, for any even $k,\ell\geq 4$ not both equal to $4$ there is an orientably-regular map for infinitely many values of $q$, cf. \cite{CPS2}. Our next result shows that this fails to hold in the case of $M(q^2)$.

\begin{prop}
If $k,\ell\equiv 0\pmod{8}$ and $k\not\equiv\ell\pmod{16}$ then there is no orientably-regular map of type $(k,\ell)$ on $M(q^2)$ for any $q$.
\end{prop}

\pr
By Theorem \ref{t-conj} orders of elements in $G{\setminus}G_0$ are $o_i=2(q-1)/{\rm gcd\{q-1,i\}}$ and $o'_i=2(q+1)/{\rm gcd\{q+1,i\}}$ for odd $i$ such that $1\le i\le (q-1)/2$ and $1\le i\le (q+1)/2$, respectively. Note that if $o_i\equiv 0\mod 8$ then $o'_i\equiv 4\mod 8$ and vice versa. Further, if $o_i\equiv 8\mod 16$ then $q-1\equiv 4\mod 8$ since $i$ is odd, and if $o_i\equiv 0\mod 16$ then $q-1\equiv 0\mod 8$. It follows that for a given $q$ we cannot have a non-singular generating pair of orders $o_i\equiv 0\mod 16$ and $o_j\equiv 8\mod 16$. The argument for orders of the form $o'_i$ is similar.
\ebox

Besides reflexibility, another frequently studied property of orientably-regular maps is self-duality.
In general, an oriented map is \emph{positively self-dual} if it is isomorphic to its dual with the same orientation, and \emph{negatively self-dual} if it is isomorphic to its oppositely oriented dual map. In terms of orientably-regular maps represented by triples $(G,x,y)$, positive and negative self-duality is equivalent to the existence of an (involutory) automorphism of $G$ sending the ordered pair $(x,y)$ onto $(y,x)$ and $(y^{-1},x^{-1})$, respectively. By the same arguments as in the second paragraph of Section \ref{sec:ref} one concludes that for our group $G=M(q^2)$, an orientably-regular map defined by a generating pair $[A,1],[B,1]$ will be positively self-dual if and only if there exists an involution $[C,i]\in\ovl{G}$ conjugating the two generators, and the map will be negatively self-dual if there is such an involution conjugating $[A,1]$ to $[B,1]^{-1}$.
\smallskip

Setting up the corresponding matrix equations for such conjugations, however, lead to enormously complicated formulae from which we were not able to extract `nice' closed formulae. Clearly in a self-dual map we have $k=\ell$ so that the orders of $[A,1]$ and $[B,1]$ must be equal. We used GAP\cite{GAP} to construct all the regular maps on $M(q^2)$, for small values of $q$, with the generators of equal order, and then tested for self-duality by determining if a conjugating element $[C,i]$ as above exists. The results of this computation are given in Table~\ref{tab:selfdual}, showing the numbers $\#(k{=}\ell)$ of maps that have generators of equal orders, those which are positively or negatively self-dual, and those which are both.

\begin{table}[h]
\begin{tabular}{|c|cccc|cccc|}
\hline
&\multicolumn{4}{|c|}{$B=\mathrm{dia}(\lambda,1)$}&\multicolumn{4}{|c|}{$B=\mathrm{off}(\lambda,1)$}\\
$q$ & $\#(k{=}\ell)$ & $+$ self-dual & $-$ self-dual & both & $\#(k{=}\ell)$ & $+$ self-dual & $-$ self-dual & both\\
\hline
3 & 0 & 0 & 0 & 0 & 3 & 3 & 3 & 3\\
5 & 15 & 15 & 5 & 5 & 10 & 10 & 6 & 6\\
7 & 28 & 28 & 8 & 8 & 78 & 42 & 14 & 14\\
9 & 95 & 45 & 9 & 9 & 68 & 36 & 10 & 10\\
11 & 276 & 132 & 24 & 24 & 265 & 165 & 33 & 33\\
13 & 469 & 273 & 39 & 39 & 666 & 234 & 42 & 42\\
17 & 2556 & 612 & 68 & 68 & 1312 & 544 & 72 & 72\\
19 & 1960 & 760 & 80 & 80 & 2799 & 855 & 95 & 95\\
\hline
\end{tabular}
\caption{Numbers of self-dual maps on $M(q^2)$}
\label{tab:selfdual}
\end{table}

Note that the computational evidence suggests that a negatively self-dual map on $M(q^2)$ is also positively self-dual.
\bigskip
\bigskip

\noindent{\bf Acknowledgement.}~ The second and the third authors acknowledge support from the APVV 0136/12 and APVV-015-0220 research grants, as well as from the VEGA 1/0007/14 and 1/0026/16 research grants.

\end{document}